\documentclass[11 pt,a4paper,twoside]{article}
\usepackage{amsmath,amssymb,amsthm,epsfig,graphics}

\renewcommand{\AA}{\mathbb{A}}

\newcommand{\NN}{\mathbb{N}}

\newcommand{\RR}{\mathbb{R}}
\renewcommand{\SS}{\mathbb{S}}
\newcommand{\TT}{\mathbb{T}}

\newcommand{\ZZ}{\mathbb{Z}}

    \def\cU{{\cal U}}
    \def\cV{{\cal V}}
    
\def\cF{{\cal F}}

\def\ra{\rightarrow}

\renewcommand{\phi}{\varphi}
\renewcommand{\epsilon}{\varepsilon}

\newtheorem{theo}{Theorem}[section]
\newtheorem{prop}[theo]{Proposition}
\newtheorem{coro}[theo]{Corollary}
\newtheorem*{conj*}{Conjecture}

\newtheorem{lemmanota}[theo]{Lemma and notation}

\newtheorem{rema}[theo]{Remark}

\newtheorem{ques}[theo]{Question}
\newtheorem*{prob*}{Problem}

\renewcommand{\t}{\widetilde}
\def\wt{\widetilde}

\newtheorem*{theo*}{Theorem}

\def\adhe{\mathrm{Cl}}

\def\rot{\mathrm{Rot}}
\def\leb{\mathrm{Leb}}
\def\fix{\mathrm{Fix}}

\newcommand{\homeo}{\operatorname{Homeo}}

\title{Pseudo-rotations of the open annulus}
\author{F. B\'eguin\footnote{Laboratoire de math\'ematiques,
    Universit\'e Paris Sud,
91405 Orsay Cedex, France.}, S. Crovisier\footnote{CNRS - Laboratoire Analyse, G\'eom\'etrie et Applications, UMR 7539, Institut Galil\'ee,
Universit\'e Paris 13, Avenue J.-B. Cl\'ement, 93430 Villetaneuse, France.}
and F. Le Roux\footnote{Laboratoire de math\'ematiques, Universit\'e Paris Sud, 
91405 Orsay Cedex, France.}}

\begin{document}

\sloppy

\maketitle

\begin{abstract}
In  this  paper,  we   study  pseudo-rotations  of  the  open  annulus,
\emph{i.e.}  conservative  homeomorphisms  of  the open  annulus  whose
rotation set is  reduced to a single irrational  number (the angle of
the  pseudo-rotation).   We  prove   in  particular  that,  for  every
pseudo-rotation  $h$ of  angle  $\rho$, the  rigid  rotation of  angle
$\rho$ is in the closure of  the conjugacy class of $h$. We also prove
that  pseudo-rotations are not  persistent in  $C^r$ topology  for any
$r\geq 0$.
\end{abstract}

\section*{Introduction}
\label{s.intro}


\subsection{Some motivations}

The    concept    of     rotation    number    was    introduced    by
H. Poincar\'e~\cite{poincare}  to compare the  dynamics of orientation
preserving  homeomorphisms of  the  circle to  the  dynamics of  rigid
rotations.   To  any   orientation-preserving  homeomorphism   $h$  is
associated a unique rotation number $\rho(h)$, measuring in some sense
the average speed of rotation of the orbits of $h$ around the circle. 
In  the case  where $\rho(h)$  is rational,  the dynamics  of  $h$ may
degenerate  dramatically:  $h$ may  present  only  one periodic  orbit
(whereas, for the rigid rotation $R_{\rho(h)}$, all the orbits are
periodic). On the contrary, in the case where $\rho(h)$ is irrational,
$h$ is always semi-conjugate  to the rigid rotation $R_{\rho(h)}$, and
the closure  of the  conjugacy class of  $h$ always coincides  with the
closure of the conjugacy class of the rotation $R_{\rho(h)}$.

The notion of rotation  number was generalized by Misiurewicz, Ziemian,
and Franks in order to  describe the dynamics of homeomorphisms of
the closed annulus and of the two-torus (see e.g. \cite{MisZie}). 
More recently, it was used by P.~Le~Calvez in order to describe
the dynamics of conservative homeomorphisms of the open annulus. Given
a homeomorphism $h$ of  the (closed or open) annulus  isotopic to
the identity, one can 
define the \emph{rotation set} of $h$,  which is in some sense the set
of all the possible asymptotic speeds  of rotation of the orbits of $h$
around  the annulus.  This is  a subset  of $\RR$,  defined up  to the
addition of  an integer. In  general, the rotation  set of $h$  is not
reduced to a single point, and the dynamics of $h$ is much richer than
the dynamics  of a  single rotation.  However, one  can address  the following
problem:

\begin{prob*}
Consider a  homeomorphism $h$ of  the annulus, such that  the rotation
set  of   $h$  is  reduced  to   a  single  number   $\rho$  which  is
irrational (such   an   homeomorphism   will  be   called   a
  \emph{pseudo-rotation of angle $\rho$}).  To what extend does the 
dynamics of $h$ looks like the rigid rotation with angle $\rho$~?  
\end{prob*}

In  the case  of the  closed annulus  $\SS^1\times [-1,1]$,  the above
problem  has   been  studied  in~\cite{beguin-crovisier-leroux-patou},
starting     from    a    generalization     of    a     theorem    of
J.~Kwapisz~\cite{kwapisz}. We would like to deal here with the case of
the open annulus $\SS^1\times\RR$.

Results on homeomorphisms of the  open annulus are usually much harder
to prove than their analogs  on the compact annulus. However, the open
annulus setting has a particular interest: 
it is related to the conservative dynamics on the two-sphere.  Indeed,
any  orientation-preserving  conservative  homeomorphism $h$  of  the
two-sphere $\SS^2$ has at least two distinct fixed points $N$ and $S$;
removing these two points, one gets a homeomorphism of
the    open    annulus    $\SS^2\setminus\{N,S\}\simeq\SS^1\times\RR$.
Moreover,  the rotation  set of  this  homeomorphism is  reduced to  a
single  irrational number if  and only  if $h$  has no  other periodic
points  than $N$  and  $S$ (see  proposition~\ref{p.no-periodic-orbit}).
This is  the reason the above-mentionned problem is  connected to the
following  conjecture  of  G.~Birkhoff (see~\cite[page  712]{birkhoff}
and~\cite{herman2}).
\begin{conj*}[Birkhoff's sphere conjecture]
 Let  $h$  be  an  orientation preserving  real-analytic  conservative
 diffeomorphism  of  the  two-sphere, and having  only  two  periodic
 (necessarily fixed) points. Then, $h$ is conjugate to a rigid rotation.  
\end{conj*} 

This  conjecture is  still open. An  example of M. Handel, improved by
M.~Herman, shows that  the real-analyticity  assumption is  necessary:
there  exists a 
$C^\infty$ diffeomorphism  of the two-sphere, having  only two periodic
(fixed)   points,  that  is   not  conjugate   to  a   rigid  rotation
(\cite{handel,herman1}).  Note  that,   in
Handel-Herman construction, the rotation number of the diffeomorphism 
is  necessarily a Liouville number.  On the  contrary,  in the  case were  the
rotation  number is assumed  to be  diophantian, some  partial results
towards the conjecture,  based on KAM theory and working for
$C^\infty$  diffeomorphisms, were  proposed by Herman
and  written in~\cite{fayad-krikorian-vivier}. 
Our results, far from proving the conjecture, give some
kind of qualitative and topological motivation for it.  

\subsection{The line translation theorem}
Let  us  denote  by  $\AA=\SS^1\times\RR$  the open  annulus.  We  can
identify $\AA$ with  the sphere $\SS^2$ minus two  points $N$ and $S$.
We call \emph{Lebesgue probability measure} on $\AA$ the measure induced by
the Lebesgue measure  on $\SS^2$. We call \emph{essential
topological  line in  $\AA$}  every simple  curve,  parametrized by  $\RR$,
properly embedded in  $\AA$, joining one of the ends  of $\AA$ to the
other. We  recall that a \emph{Farey  interval} is an  interval of the
form   $]\frac{p}{q},\frac{p'}{q'}[$    with   $p,q,p',q'\in\ZZ$   and
$qp'-pq'=1$.  Here is our main result.

\begin{theo}[Line translation theorem]
\label{th.arc-translation}
Let $h\colon \AA\to\AA$ be a homeomorphism of the open annulus which is
isotopic   to    the   identity,   which    preserves   the   Lebesgue
measure.   Assume that the  closure of  the rotation  set of  some
lift $\t h:\RR^2\rightarrow\RR^2$ of $h$ is 
contained in a Farey interval $]\frac{p}{q},\frac{p'}{q'}[$.  

Then, there  exists an  essential topological line  $\gamma$ of  $\AA$ such
that  the topological lines $\gamma,h(\gamma),\dots,h^{q+q'-1}(\gamma)$
are pairwise  disjoint. Moreover, the  cyclic order of these topological lines is
the  same as  the cyclic  order of  the $q+q'-1$  first iterates  of a
vertical line $\{\theta\}\times \RR$  under the rigid rotation with
angle $\rho$, for any $\rho\in]\frac{p}{q},\frac{p'}{q'}[$.  
\end{theo}

Very  roughly speaking, theorem~\ref{th.arc-translation}  asserts that,
if the  the rotation set  of a homeomorphism  $h:\AA\rightarrow\AA$ is
included in  a Farey interval  $]\frac{p}{q},\frac{p'}{q'}[$, then the
dynamics  of $h$  is similar  to those  of a  rigid rotation  of angle
$\rho\in]\frac{p}{q},\frac{p'}{q'}[$, provided that  one does not wait
for more than $q+q'-1$ iterates.

Although  the  statement  of theorem~\ref{th.arc-translation}  is  the
natural generalization of the arc translation theorem
of~\cite{beguin-crovisier-leroux-patou},  the proofs  of  these  two
results are  completely different. Indeed, most of  the arguments used
in~\cite{beguin-crovisier-leroux-patou}  are specific  to  the compact
annulus~; here, we will have  to use some technics coming from Brouwer
theory,  that are  typical from topological  dynamics on  non-compact
surfaces.

The line translation theorem implies the following useful corollary: \emph{if the rotation set of $h$ is bounded, then $h$ is conjugate to a homeomorphism whose displacement function is bounded} (see proposition~\ref{p.integration} below). This corollary plays a key role in the proof of the perturbation theorem~\ref{t.perturbation-pseudo-rotation} below.

\subsection{Results on pseudo-rotations}

We call  \emph{pseudo-rotation} of the open  annulus any homeomorphism
which  is  isotopic to  the  identity,  which  preserves the  Lebesgue
measure,  and  whose  rotation  set  is reduced  to  a  single  number
$\alpha$.   This number  $\alpha$ (defined  up to  the addition  of an
integer)  is  called the  {\em  angle}  of  the pseudo-rotation.   The
following   proposition   provides   an  alternative   definition   of
pseudo-rotations with irrational angles: 

\begin{prop}[Characterization of pseudo-rotations] 
\label{p.no-periodic-orbit}
Let $h$ be a homeomorphism  of the open  annulus $\AA$, isotopic to
the identity and  preserving the  Lebesgue probability measure. then
$h$   is a  pseudo-rotation with  irrational angle if and only if it
does not have any periodic orbit.  
\end{prop}

This result does not seem to appear in the literature. It can be seen
as a straightforward application of a generalisation of
Poincar\'e-Birkhoff theorem by J. Franks, together with an ergodic
theoretical argument of P. Le Calvez. 
We will provide a proof in section~\ref{ss.existence-periodic}.

As    an    immediate    corollary    of    the    line    translation
theorem~\ref{th.arc-translation}, we get: 

\begin{coro}[Line translation theorem for pseudo-rotations]
\label{c.line-translation}
Let  $h:\AA\rightarrow\AA$ be a  pseudo-rotation of  irrational angle
$\rho$.  Then,  for  every  $n\in\NN\setminus\{0\}$, there  exists  an
essential  topological line  $\gamma$  in $\AA$,  such  that the
topological  lines   $\gamma,h(\gamma),\dots,h^n(\gamma)$   are
pairwise 
disjoint.      The      cyclic       order      of      the      lines
$\gamma,h(\gamma),\dots,h^n(\gamma)$ is the same as the cyclic order of
the  $n$  first iterates  of  a  vertical  line under  the  rigid
rotation of angle $\rho$. 
\end{coro}

Corollary~\ref{c.line-translation} can be seen as an analogue of the
following well-known property  for the dynamics on the  circle: if $h$
is  an   orientation-preserving  homeomorphism  of   the  circle  with
irrational  rotation number  $\rho$,  then the  cyclic  order of  the
points of  any orbit of  $h$ is  the same as  the cyclic order  of the
points of any orbit of the rigid rotation with angle $\rho$. 
However,   note  that,   in   corollary~\ref{c.line-translation},  the
essential  simple  line  $\gamma$  does  depend on  the  integer  $n$.
Indeed, one can construct a pseudo-rotation $h:\AA\rightarrow\AA$ with
irrational  angle such  that no  essential  topological line  in $\AA$  is
disjoint  from  all  its  iterates  under $h$  (see  the  examples  of
Handel~\cite{handel} and Herman~\cite{herman2}).  

Using corollary~\ref{c.line-translation}, one can prove the following: 

\begin{theo}[Closure of the conjugacy class of a pseudo-rotation]
\label{t.closure-conjugacy-class}
Let $h$ be a pseudo-rotation of the open annulus with irrational angle
$\rho$.  The rigid  rotation of angle $\rho$ is in  the closure (for
the compact-open topology) of the
conjugacy class\footnote{Here, the  conjugating homeomorphisms
are  not assumed to  be conservative.} of $h$. 
\end{theo}

In other words,  for every pseudo-rotation $h$ of  angle $\rho$, there
are  conjugates   of  $h$  which   are  arbitrarily  close   (for  the
compact-open topology) to a rigid rotation. We do not know if the same
result holds if one allows only conservative conjugacies.  We also do
not know if  any pseudo-rotation of angle $\rho$ is  in the closure of
the rigid rotation of angle $\rho$.

Corollary~\ref{c.line-translation} and
theorem~\ref{t.closure-conjugacy-class}  show   some  common  features
between the dynamics of any pseudo-rotation with irrational angle and 
the  dynamics  on a  rigid  rotation. Nevertheless, there  are
examples of  pseudo-rotations whose  dynamics is quite  different from
those  of  a rotation.   Indeed,  using  techniques developed  by
D.~Anosov     and    A.~Katok    (see~\cite{anosov-katok,fathi-herman,
  fayad-katok,fayad-saprykina}), one can construct $C^\infty$
pseudo-rotations for which the Lebesgue probability measure is
ergodic; in particular, such pseudo-rotations are not semi-conjugate
to a rigid rotation.

We end  this discussion on pseudo-rotations by  noting that irrational
pseudo-rotation are  not robust  under perturbations: for  each $r\geq
0$, the set of irrational  pseudo-rotations is meagre in the space of
$C^r$     conservative     diffeomorphisms     isotopic     to     the
identity (see Corollary~\ref{c.meager}). This will be a consequence of the
following perturbation result, where the perturbation is chosen
\emph{a priori}, and does not depend on the map one wants to perturb. 

\begin{theo}[Perturbation of pseudo-rotation]
\label{t.perturbation-pseudo-rotation} 
For every  homeomorphism $h\colon \AA\to\AA$ isotopic  to the identity
and preserving the Lebesgue probability measure, there exists a rigid
rotation $R$ of arbitrarily small angle such that $h\circ R$ has a
periodic orbit.  
\end{theo}

Theorem~\ref{t.perturbation-pseudo-rotation}  answers  a  question  of
  J.~Franks, who also proved that the same statement holds
 in  the compact annulus (see~\cite{franks2}   pages    18--19). 
 To cope with the lack of compactness, we have to use
 the line translation theorem and some continuity results of P. Le Calvez.
     Note   that   the   analogue   of
  theorem~\ref{t.perturbation-pseudo-rotation} for 
non-conservative homeomorphism of the  annulus was shown to be false
  G.~Hall  and  M.~Turpin  \cite{hall-turpin}. Moreover, it is not
  known (see~\cite{herman2}) if, for $r\geq 
  2$,  the space  of $C^r$  diffeomorphisms of  the two-torus  (in the
  non-conservative case) or of compact manifolds with dimension larger
  or equal to $3$ (in the conservative and non-conservative cases) has
  a dense subset of diffeomorphisms that present a periodic orbit. 

In a forthcoming paper, we shall prove that any irrational pseudo-rotation $h$
possesses a circle compactification  in the following sense : there exists a
homeomorphism $\hat h$ of the compact annulus $\SS^1 \times [0,1]$ 
 whose restriction to the open annulus $\SS^1 \times ]0,1[$
 is conjugate to $h$. In other words, if we see $h$ as a homeomorphism of the sphere fixing the North and South poles, one can construct a blow-up of $h$ at each fixed point.

\paragraph{Acknowledgements.}  We  would like to  thank P.\,Le\,Calvez
for many useful discussions.

\section{Preliminaries (I) : rotation numbers}
\label{s.nbrotation}


\subsection{The open annulus}

We denote by $\AA=\TT^1\times \RR$ the infinite annulus and by $\t 
\AA=\RR\times\RR$  its  universal  cover.   We  denote  by  $\pi$  the
canonical projection  of $\t \AA$ onto  $\AA$.  We denote  by $p_1$
the projection defined on $\AA$ or $\widetilde\AA$ by $p_1(x,y)=x$.
We  denote  by $T\colon  \t   \AA\to  \t   \AA$ the  translation
defined  by  $T(x,y)=(x+1,y)$. Note  that  the  annulus  $\AA$ is  the
quotient space $\t  \AA/T$.  We will sometimes consider the annulus
$\AA_q=\RR^2/T^q$ for some $q\geq 2$.  

By the two points compactification, one can identify the annulus $\AA$
to the punctured sphere $\SS^2\setminus\{N,S\}$, where $N$ and $S$ are
two  distinct points  of  $\SS^2$.  The  Lebesgue  measure on  $\SS^2$
induces on $\AA$ a probability measure  on $\AA$ that we call the {\em
Lebesgue probability measure of $\AA$} and denote by $\leb$. 


The  set  of  the  homeomorphisms   of  the  annulus  (resp.   of  the
two-sphere)  that   are  isotopic  to  the  identity   is  denoted  by
$\homeo^+(\AA)$ (resp by  $\homeo^+(\SS^2)$).  We will mostly consider
the   subsets  $\homeo^+_\leb(\AA)$   and   $\homeo^+_\leb(\SS^2)$  of
$\homeo^+(\AA)$ and $\homeo^+(\SS^2)$ made of the homeomorphisms which
preserve the Lebesgue probability measure.  



\subsection{Rotation numbers of  points and measures, rotation set of a
homeomorphism} 
\label{ss.definition}

Consider  a  homeomorphism  $h\in\homeo^+(\AA)$,  and a  lift  $\t 
h\colon \t  \AA\to \t  \AA$  of $h$. Since $\AA$ is not compact,
the definitions of the rotation number of a point under $\t  h$, of
the rotation set of $\t  h$,  \emph{etc.} cannot be as simple as in
the  case of  the  closed  annulus.  We  follow  here the  definitions
proposed by Le~Calvez in~\cite{lecalvez}.  


Let us consider a recurrent point $z\in\AA$ of $h$. We say that 
the \emph{rotation  number of $z$ under $\t  h$}  is well-defined and equal
to $\rho(z,\t  h)\in  \RR \cup \{\pm \infty\}$ if, for every lift $\t   z$ of $z$ and
for any subsequence $(h^{n_k})_{k\geq  0}$ of $(h^n)_{n\geq 0}$ and of
$(h^n)_{n\leq 0}$ such that $h^{n_k}(z)$ converges to $z$, we have 
$$
\frac{p_1\circ                 \t                  h^{n_k}(\t 
  z)}{n_k}\longrightarrow \rho(z,\t  h).
  $$

The \emph{rotation set} $\rot(\t  h)$ of $\t  h$ is the set
 of all rotation numbers of recurrent points of $\t 
h$. As it is discussed in~\cite{lecalvez}, we consider only recurrent
points  in   order  to  get   a  definition  which  is   invariant  by
conjugacy. Note that the rotation set may be empty.


Now, consider  a probability measure  $m$ on $\AA$ which  is invariant
under  $h$.   Note that  $m$-almost  every  point  is recurrent  under
$h$. Suppose that 
\begin{itemize}
\item  $m$-almost  every  point  $z\in  \AA$  has  a  rotation  number
  $\rho(z,\t h)$~; 
\item the function $z \mapsto \rho(z, \t  h)$ is integrable (with respect to the measure $m$). 
\end{itemize}
Then, we say  that \emph{the rotation number of the measure $m$ under
  $\t h$ is well-defined} and equal to
$$\rho(m,\t h)=\int_{\AA} \rho(z, \t h)d m.$$

In  the case  where  $m$ is  the  Lebesgue probability measure
\footnote{Or,  more
  generally, in the  case where $m$ is a  probability measure such that
  $m(U)>0$  for every open  subset $U$  of $\AA$.},  Le
  Calvez  found  a  nice condition  implying  that  the
  rotation number of $m$ is well-defined. First note that, if $z$ is a
  fixed  point of  $h$,  then the  rotation  number of  $z$ is  always
  well-defined and is an integer. Consider the set $\rot_{\fix}(\t  h)$ of
  the rotation numbers  of all the fixed points of  $h$. Then, one has
  the following result.

\begin{theo}[P.~Le~Calvez, existence of the mean rotation number]
\label{th.mean-rotation-number}
Suppose that $h$ preserves the Lebesgue probability measure, and that the
set $\rot_{\fix}(\t  h)$ is bounded.  Then
Lebesgue almost every point $\t x$ has a rotation number,
and the rotation set of $\t h$ is bounded.
In particular,  the rotation number $\rho(\leb,\t
h)$ of the Lebesgue probability measure under $\t  h$ is well-defined.  
\end{theo}

The rotation set, the rotation numbers of the points, and the rotation
numbers of the measures satisfy the following elementary properties.
\footnote{For item 3, note that a point which is recurrent for $h$ is
  also recurrent for $h^q$ for any $q$.} 

\begin{prop}\label{p.usual-properties}~

\begin{enumerate}
\item The rotation set, the rotation number of a point, and the
rotation number of a measure are invariant by conjugacy in
$\homeo^+_\leb(\AA)$. 
\item The rotation set of $T^k\circ\t  h$ is obtained by 
translating  by $k$ the  rotation set  of $\t  h$. Similarly,  for the
rotation number of a point, or the rotation number of an invariant
measure.
\item   The  rotation   set  of   $\t  h^q$   is   $q  \mathrm{Rot}(\t
  h)$.   Similarly for the  rotation number  of a  point, and  for the
  rotation number of an invariant measure. 
\end{enumerate}
\end{prop}

\subsection{The morphism property}
\label{ss.morphisms}

The \emph{horizontal displacement of $\t  h$} is  the  function
$r\colon  \AA\to\RR$  defined  as 
follows: given $z\in\AA$,  we choose a lift $\t  z$  of $z$, and we
set  $r(z)=p_1(\t  h(\t   z))-p_1(\t  z)$.  Note  that $r(z)$
does depend on the choice of $\t  z$.  If $m$ is an $h$-invariant
probability measure, and if  $r$ is $m$-integrable, Birkhoff's ergodic
theorem implies that  $m$ has a rotation number equal  to $\int r dm$.
This shows that the rotation  number of the Lebesgue probability
measure satisfies some morphism property. 

\begin{prop}
\label{p.morphism-property}
Let $h$, $g$ be two homeomorphisms of $\AA$ that are isotopic to the
identity and preserve the Lebesgue probability measure.  Let $\t  h$,
$\t  g$, $\t  h\circ \t  g$ be some lifts to $\t  \AA$ of
$h$, $g$ and $h\circ g$.

If  the  horizontal  displacement  of  $h$, $g$  and  $h\circ  g$  are
integrable for the Lebesgue probability measure, then  
$$\rho(\leb,\t  h\circ\t  g)=\rho(\leb,\t  h)+\rho(\leb,\t  g).$$
\end{prop}

In general, the horizontal displacement of a homeomorphism is
not integrable. Moreover,  one should note that the property
of the horizontal displacement being $\leb$-integrable is not invariant
by conjugacy.  We do not know if 
proposition~\ref{p.morphism-property}    is     true    without    the
integrability assumptions (see the precise question and the results in
paragraph~\ref{s.integrabilite}).

\section{Preliminaries (II) : Brouwer theory}
\label{s.Brouwer}

Every annulus homeomorphism
$h$  lifts to a homeomorphism $\wt h$ of the plane. Thus results about
the existence of fixed points can be obtained by considering
\emph{Brouwer  homeomorphisms}, which are  the orientation-preserving
fixed point free  homeomorphisms of  the plane $\RR^2$.
In  this section, we  briefly recall some of  the main
results of the theory of Brouwer homeomorphisms.

\subsection{Brouwer lines and Brouwer theorem} 
\label{ss.Brouwer-lines}
A \emph{topological line} in the plane is the image $\Gamma$ of a
proper continuous embedding from $\RR$ to $\RR^2$ (equivalently, using
Schoenflies theorem, it is the image of a Euclidean line under a
homeomorphism of the plane). 
Given a Brouwer homeomorphism $H$,  a \emph{Brouwer line} for $H$ is a
 topological line $\Gamma$, disjoint from its image $H(\Gamma)$, 
and such  that $\Gamma$
separates $H(\Gamma)$ from $H^{-1}(\Gamma)$.
We will say that $\Gamma$ is an \emph{oriented Brouwer line} if it is endowed with
the orientation such that $H(\Gamma)$ is on the right of $\Gamma$ (and thus $H^{-1}(\Gamma)$ is on the left of $\Gamma$).
Then for every $k\in\ZZ$, we can endow the line $H^k(\Gamma)$
with the image by $H^k$ of the orientation of $\Gamma$.
Since $H^k$ preserves the orientation, the line $H^{k+1}(\Gamma)$
  is on  the right of  $H^k(\Gamma)$, and the line $H^{k-1}(\Gamma)$ is  on the
  left of  $H^k(\Gamma)$.  By induction, we see  that $H^q(\Gamma)$ is
  on the right  of $H^p(\Gamma)$ if and only  if $q>p$. In particular,
  the lines $(H^k(\Gamma))_{k\in\ZZ}$ are pairwise disjoint. 

Now  let $U$ be the open region of $\RR^2$ situated between the
  lines $\Gamma$ and $H(\Gamma)$, and $\adhe(U)=\Gamma \cup U \cup H(\Gamma)$.
  The  sets $(H^k(U))_{k\in\ZZ}$
  are   pairwise  disjoint.   As a consequence,
  the   restriction  of   $H$  to the open set
  $O=\bigcup_{k\in\ZZ}H^k(\adhe(U))$ is conjugate to a translation.  In particular, if 
  the iterates of $\adhe(U)$ cover the whole plane, then $H$ itself is conjugate to a translation.

The  main  result  of  Brouwer  theory is  the  \emph{plane
translation  theorem}: \emph{every point  of $\RR^2$ lies on  a Brouwer
line   for   $H$   (see for example~\cite{guillou1})}.   

\subsection{Guillou-Sauzet-Le Calvez theorem}

In the case  where the Brouwer  homeomorphism $H$ is a
lift of a  homeomorphism of the annulus $\AA$, one  would like to have
an   ``equivariant  version''  of   the  plane   translation  theorem,
\emph{i. e.} one  would like to find  some Brouwer lines  for $H$ which
project as ``nice'' curves in  the annulus $\AA$.  This is the purpose
of a result of L. Guillou (see~\cite{guillou2}), which was improved by
A. Sauzet in his PhD thesis (see~\cite{sauzet}). We give below a foliated
version of Guillou-Sauzet's result which relies on a recent and
powerful theorem of P. Le Calvez (see~\cite{lecalvez2}). For sake of
simplicity, we restrict ourselves to the case of homeomorphisms
without wandering points. 
Recall that  an \emph{essential topological line} is the image of the
line $\{0\} \times \RR$ under a homeomorphism of the annulus that is
isotopic to the identity. 

\begin{theo}[L. Guillou, A. Sauzet, P. Le Calvez]
\label{t.Guillou-Sauzet}
Let   $h:\AA\rightarrow\AA$  be  a homeomorphism isotopic  to  the
identity. Assume that: 
\begin{itemize}
\item $\wt h:\RR^2\rightarrow\RR^2$
is a fixed  point  free lift of $h$; 
\item the homeomorphism $h$ does not have any wandering point
(\emph{i.e.} every open set must meet some of its iterates under $h$).
\end{itemize} 
Then there exists an oriented foliation $\cF$ of
 the annulus $\AA$ such that each oriented leaf of 
$\cF$ is an   essential topological line
  which lifts in $\RR^2$ to an oriented Brouwer line for $\wt h$. 
\end{theo}
Note that any foliation of the annulus by essential topological lines is homeomorphic to the trivial foliation by vertical lines.


\begin{proof}[Proof of theorem~\ref{t.Guillou-Sauzet}]
Let  $h$  be  a homeomorphism  of  the  annulus  $\AA$, and  let  $\wt
h:\RR^2\rightarrow\RR^2$ be a fixed point  free lift of $h$. Le Calvez
has  proved that there exists a $C^0$ oriented foliation $\cF$ of the
annulus $\AA$,  which lifts  as an oriented foliation  $\wt\cF$ of  $\RR^2$ such
that  every oriented leaf   of  $\wt\cF$  is  an oriented  Brouwer   line  $\Gamma$ for
$\wt  h$, with $\wt h(\Gamma)$ on the right of $\Gamma$ 
(see~\cite{lecalvez2}). Now we see
the annulus  $\AA$ as the  sphere minus the two  points $N,S$, and  we see
$\cF$  as  a foliation  of  $\SS^2$  with  two singularities  $N$  and
$S$.

Suppose that $\cF$ has a leaf $\gamma$ which is homeomorphic to a
circle. Since it lifts to a topological line $\Gamma$ in the universal
covering of $\AA$, this leaf must separate $N$ and $S$.  Since
$\Gamma$ is a Brouwer line, the leaf $\gamma$ is disjoint from its
image, and the open annular region $U$ between $\gamma$ and
$h(\gamma)$ is disjoint from its iterates under $h$, which contradicts
the second assumption of the theorem. 

Similarly, we see that $\cF$ does not admit a leaf which is closed in
$\AA$ and whose endpoints in $\SS^2$ are both equal to $N$, or both
equal to $S$. 
Nor does $\cF$ admits  any cycle of oriented leaves $\gamma_{1},
\gamma_{2}$ that are closed in $\AA$ and goes respectively from $N$ to
$S$ and from $S$ to $N$. Now Poincar\'e-Bendixson theory tells us that
all the leaves of $\cF$ are closed in $\AA$, and either they all go
from $N$ to $S$, or they all go from $S$ to $N$. 
\end{proof}

\begin{rema}
In most situations, we will not need the whole foliation provided
by theorem~\ref{t.Guillou-Sauzet} but only one leaf of this foliation.
\end{rema}

\subsection{Application to the existence of periodic orbits}
\label{ss.existence-periodic} 

In this section, we use some Guillou-Sauzet-Le Calvez theorem to prove
classical results about the existence of 
periodic orbits. In particular, we provide the characterisation of
irrational pseudo-rotations announced in the introduction, namely that
an annulus homeomorphism does not have any periodic orbit if and only
if its rotation set is reduced to a single irrational number
(proposition~\ref{p.no-periodic-orbit}). 


\begin{theo}[Franks~\cite{franks2}, Le Calvez~\cite{lecalvez}]
\label{t.PBFLC}
Let $h \in \homeo^+_\leb(\AA)$, and let $\t h$ be a lift of $h$.
Suppose that $\t h$ does not have any fixed point. Then the rotation
set $\rot(\t h)$ is either contained in $[-\infty,0]$ or in $[0,
+\infty]$. Furthermore, Lebesgue almost every recurrent point has a
non zero rotation number. 
\end{theo}

We do not know if the statement can be improved by proving that the
rotation set does not contain zero. 

\begin{proof} 
Let $h \in \homeo^+_\leb(\AA)$, and let $\t h$ be a lift of $h$ that
has no fixed point. Let $\t \cF$ be the lift to $\RR^2$ of the oriented foliation $\cF$ provided by theorem~\ref{t.Guillou-Sauzet}.
Either all the leaves of $\cF$ are oriented from $S$ to $N$, or they are all oriented from $N$ to $S$. In the remainder, we assume that we are in the first situation. We will prove that the rotation set of $\wt h$ is contained in
$[0,+\infty]$ and that Lebesgue almost every point has a positive  
rotation number. 

Let $\Gamma, \Gamma'$ be  lifts in $\RR^2$ of essential topological lines (oriented from $S$ to $N$).
We denote by
$L(\Gamma)$ the connected component of $\RR^2\setminus\Gamma$ 
on the left of $\Gamma$, and by $R(\Gamma)$ the connected component of
$\RR^2\setminus\Gamma$  {on the right} of $\Gamma$. We will write $\Gamma < \Gamma'$
if $\Gamma'$ is included in $R(\Gamma)$.

 Observe
that, due to the orientations, for every $\Gamma\in\wt\cF$, and every $p,q \geq 0$, 
$$ T^{-p}(\Gamma)   < \Gamma < T^{p}(\Gamma)  \;\;\;\;\;\;\;  \mbox{ and }\;\;\;\;\;\;\; 
  {\t h}^{-q} (\Gamma)  < \Gamma < {\t h}^{q}(\Gamma).$$


Consider a point $x\in\RR^2$ and a leaf $\Gamma$ of $\wt\cF$ such that
$x\in R(\Gamma)\cap L(T(\Gamma))$. On the one hand, for every $q\geq 0$, the point $\wt
h^q(x)$ is in $\wt h^q(R(\Gamma))\subset R(\Gamma)$. On the other
hand, for every $p>0$, the point $T^{-p}(x)$ is in 
$T^{-p}(L(T(\Gamma))=
T^{-p+1}(L(\Gamma))\subset L(\Gamma)$. This implies
that, the point $x$ cannot have a negative rotation 
number. This proves that the rotation set of $\wt h$ is included in $[0,+\infty]$.

We are left to prove that Lebesgue almost every point in $\RR^2$ has a
positive rotation number. For this purpose, we use some ergodic
theoretical arguments due to P. Le Calvez (see~\cite[page 3227]{lecalvez}).
Consider a leaf $\Gamma$ of $\wt\cF$. 
Let 
$$\wt U=\wt U_\Gamma=R(\Gamma)\cap L(\wt h(\Gamma))\cap
 L(T(\Gamma)),$$
and $U=U_\Gamma$ be the projection in $\AA$ of $\wt U$. Note that,
by definition, $\wt U$ is disjoint from its images under $\wt h$
and $T$. 
Consider the \emph{return time function} $\nu=\nu_\Gamma:U\rightarrow\NN\setminus \{0\}$,
the \emph{first return map} $\Phi=\Phi_\Gamma:U\rightarrow U$, and the
\emph{displacement function} $\tau=\tau_\Gamma:U\rightarrow\ZZ$
defined as follows: 
\begin{itemize}
\item $\nu(x)=\inf\{n > 0 \mid h^n(x)\in U\}$;
\item $\Phi(x)=h^{\nu(x)}(x)$;
\item  $ \tau(x) $  is the unique integer such that $ \wt
  h^{\nu(x)}(\wt x)\in T^{\tau(x)}(\wt U)$, where $\wt x$ is the (unique) lift of $x$ in $\wt U$.
\end{itemize}  
By classical arguments, the function $\nu$ is
integrable. Hence, by Birkhoff ergodic theorem, the quantity    
$$\nu^*(x)=\lim_{n\rightarrow
  +\infty}\frac{1}{n}\sum_{k=0}^{n-1}\nu(\Phi^k(x))$$   
exists, is finite and positive for Lebesgue almost every $x$ in
  $U$. We claim that $\tau(x)$ is a positive integer for
  every $x\in U$:  indeed, for every $\wt x\in\wt U$,
  the point $\wt 
  h^{\nu(x)}(\wt x)$ is in $\wt h^{\nu(x)}(R(\Gamma))$,
  which is included in $R(\wt h(\Gamma))$,
  and, for every $p\geq 0$, the set $T^{-p}(\wt U)$ is contained
  in $L(\wt h(\Gamma))$. Hence, by Birkhoff
  ergodic theorem for positive functions, the quantity
  $$\tau^*(x)=\lim_{n\rightarrow 
  +\infty}\frac{1}{n}\sum_{k=0}^{n-1}\tau(\Phi^k(x))$$  
exists and is  greater than or equal to $1$ (maybe equal to $+\infty$) for Lebesgue almost
  every $x$ in $U$.   
Since $U$ is open, the recurrent points of $h$ in $U$
  are exactly the recurrent points of $\Phi$. Hence, the
  rotation number of Lebesgue
  almost every point $x$ of $U$ is equal to 
$$\lim_{n\rightarrow +\infty}
\frac{\tau(x)+\dots+\tau(\Phi^{n-1}(x))} 
{\nu(x)+\dots+\nu(\Phi^{n-1}(x))}=
\frac{\tau^*(x)}{\nu^*(x)},$$
which is positive (maybe equal to $+\infty$) for Lebesgue almost every
point in $U$.
Since $\RR^2=\bigcup_{\Gamma\in\wt\cF} U_\Gamma$, and 
since $U_\Gamma$ is a non-empty open set for every $\Gamma$, this
implies that almost every point in $\RR^2$ has a non-zero rotation number.
\end{proof}


\begin{coro}
Let $h \in \homeo^+_\leb(\AA)$, and let $\t h$ be a lift of $h$.
Let $\frac{p}{q}$ be a rational number in $]\rho^-,\rho^+[$, where
$\rho^-$ and $\rho^+$ belong to the rotation set of $\t h$.
Then $\frac{p}{q}$ also belongs to the rotation set, and is the rotation number of a $q$-periodic point of $h$.
\end{coro}

\begin{proof}
Apply the previous theorem to $T^{-p} \t h ^{q}$ (using proposition~\ref{p.usual-properties}).
\end{proof}

\begin{proof}[Proof of proposition~\ref{p.no-periodic-orbit}]
Let $h \in \homeo^+_\leb(\AA)$, and let $\t h$ be a lift of $h$.
Any periodic point of $h$ has a rational rotation number, which proves
the easy part of the proposition. So assume that $h$ does not have any
periodic orbit. 
According to the previous corollary, the rotation set of $\t h$ is reduced to a single number $\alpha$.
Furthermore, the second part of
theorem~\ref{t.PBFLC} (applied to the homeomorphisms $T^{-p} \t h ^{q}$) implies that $\alpha$
 cannot be a rational number. This
completes the proof. 
\end{proof}

%

\section{Proof of the line translation theorem}
\label{s.preuvearctranslate} 

The  purpose  of  this  section   is  to  prove  the  line  translation
theorem~\ref{th.arc-translation}.  Let us explain briefly the strategy
of the proof.  In subsection~\ref{ss.conjugue-a-translation}, we prove
a preliminary result which  ensures that a homeomorphism whose
rotation set is contained in $[\epsilon,+\infty[$ for some
$\epsilon>0$ is conjugate to a  translation. In
subsection~\ref{ss.proposition-arithmetique}, we 
introduce  the  first  return  maps $\wt\phi=T^{-p}\circ\wt  h^q$  and
$\wt\psi=T^{p'}\circ\wt  h^{-q'}$, and we  state a  proposition saying
that, to prove theorem~\ref{th.arc-translation},  it is enough to find
an essential simple {line} $\gamma$ in $\AA$ and a lift of $\gamma$ which
is  disjoint from  its images  under $\wt\phi$  and  $\wt\psi$.  This
proposition is  a classical consequence of  arithmetical properties of
Farey  intervals.   Subsection~\ref{ss.proof-arc-translation} contains
the core of 
the proof  of theorem~\ref{th.arc-translation}.  The
results of subsection~\ref{ss.conjugue-a-translation} implies that the
homeomorphism  $\wt\phi$ is conjugate  to a  translation, so  that the
quotient  $\AA':=\RR^2/\wt\phi$  is homeomorphic  to  an annulus.  The
homeomorphism  $\wt\psi$  induces  an  homeomorphism  $\psi'$  of  the
annulus  $\AA'$.  So,  we  can  apply Guillou-Sauzet-Le Calvez theorem
to  the
homeomorphism  $\psi'$.  It provides  us  with a  line $\Gamma$  in
$\RR^2$, which is a Brouwer line for $\wt \psi$, 
and projects in $\AA'$ as an essential topological {line}.
 Thus $\Gamma$ is  is also a Brouwer line for 
$\wt \phi$. Then, we prove that $\Gamma$ is also
a  Brouwer line  for the  translation $T$, and that it projects to 
an essential topological {line} in our  original annulus $\AA$.

Note that the we do not know if one can strengthen the statement of theorem~\ref{th.arc-translation} by removing the word \emph{closure}. The
example described in appendix~\ref{s.exemple-a-poils} only shows that our strategy fails to prove this stronger result, since the first step of the proof (proposition~\ref{p.conjugue-a-translation} below) does not work anymore.


\subsection{Homeomorphisms with positive rotation sets} 
\label{ss.conjugue-a-translation}
The purpose of this subsection is to prove the following.
\begin{prop}
\label{p.conjugue-a-translation}
Let   $g \in \homeo^+_{\leb}(\AA)$, and $\wt g:\RR^2\rightarrow\RR^2$ be a lift
of $g$.  Assume that  the closure of  the rotation  set of $\wt  h$ is
included  in   $]0,+\infty]$.   Then  $\wt  g$  is   conjugate  to  a
translation.  
\end{prop}
Note that the  above  statement  is  sharp:  one can  construct  an  example  of
a measure-preserving homeomorphism $g:\AA\rightarrow\AA$ isotopic to the
identity, such that, for some lift $\wt g$ of $g$, the rotation set of
$\wt g$ is included in $]0,+\infty]$,  but $\wt g$ is not conjugate to
a translation (see appendix ~\ref{s.exemple-a-poils}).

\begin{proof}[Proof of proposition~\ref{p.conjugue-a-translation}]
Choose a positive  integer $k$ such that the  rotation set of $\wt g$
is included 
in   $]\frac{1}{k},+\infty]$. Consider  the homeomorphism  $\wt  g':=\wt
g^k\circ T^{-1}$, which is a lift of the homeomorphism $g'=g^k$.
 The  rotation set of $\wt g'$  is included in $]0,+\infty]$ (see
 proposition~\ref{p.usual-properties}). In particular, 
the homeomorphism $\wt g'$ is fixed  point free.
Furthermore, since $g$ preserves the Lebesgue probability measure on
$\AA$, so does $g'$, and in particular no point is wandering under the
action of $g'$.  
Thus we can apply Guillou-Sauzet theorem~\ref{t.Guillou-Sauzet}, which
provides us with an essential topological line $\gamma$ in $\AA$, such
that some lift $\Gamma$ of $\gamma$ is disjoint from its image $\t
g'(\Gamma)$.  

Using the conservative version of Schoenflies theorem (see
appendix~\ref{s.schoenflies}), we can assume that $\Gamma$ is the
vertical line $\{0\} \times  \RR$ in $\RR^2$, oriented from bottom to
top. 
The image $\t g'(\Gamma)$ is disjoint from $\Gamma$. If it was on the
left side of $\Gamma$, then the rotation set of $\t g'$ would be
contained in $[-\infty,0[$ (by the same argument as in the proof of
theorem~\ref{t.PBFLC}). Thus $\t g'(\Gamma)$  is on the right side
of $\Gamma$. Applying the covering translation $T$, we get that $\t
g^k (\Gamma)$ is on the right side of $T(\Gamma)$. By induction, for
any positive integer $n$, $\t g ^{nk} (\Gamma)$ is on the right of
$T^n(\Gamma)$. Similarly, the topological line $\t g ^{-nk} (\Gamma)$
is on the left of $T^{-n}(\Gamma)$.  
 Let $\adhe(U)$ denote the closed band delimited by $\Gamma$ and $\t
 g^k(\Gamma)$;  we get that the iterates of $\adhe(U)$ by $\t g^k$
 cover the whole plane. 
 Thus $\t g^k$ is conjugate to a translation (see
 paragraph~\ref{ss.Brouwer-lines}). 

Now it follows from a standard argument that $\t g$, having a power
conjugate to a translation, is also conjugate to a translation (the
quotient $\t \AA / \t g^k$ is an annulus, 
thus $\t \AA / \t g$ is the quotient of an annulus by a map of finite
order: this is a topological  surface whose fundamental group is
infinite cyclic, so (using the classification of surfaces) it is again
an annulus, so that $\t g$ is conjugate to a translation). 
\end{proof}

\subsection{The  ``first return maps''  $\wt\phi=T^{-p}\circ\wt h^q$
  and $\wt\psi=T^{p'}\circ\wt h^{-q'}$} 
\label{ss.proposition-arithmetique}

We consider  a homeomorphism  $h \in \homeo_{\leb}^+(\AA)$, and a  lift $\wt
h:\RR^2\rightarrow\RR^2$   of   $h$. We assume that the rotation  set
of $\wt  h$ is 
included  in   a  Farey  interval   $]\frac{p}{q},\frac{p'}{q'}[$.  We
consider   the   homeomorphisms   $\wt\phi:=T^{-p}\circ\wt  h^q$   and
$\wt\psi:=T^{p'}\circ\wt  h^{-q'}$,  sometimes  called  \emph{the
  first  return maps  associated with  $h$}. These  two homeomorphisms
play a fundamental  role in the proof of  the line translation theorem,
via the following proposition.

\begin{prop}
\label{p.arithmetics}
Let $\gamma$ be  an essential  {topological line} in  the annulus
$\AA$. Assume 
that some lift $\Gamma$ of $\gamma$ is disjoint from its images under
the   first  return   maps   $\wt\phi$  and   $\wt\psi$.  

 Then   the
$q+q'-1$ first iterates of $\gamma$  under $h$ are pairwise disjoint,
and ordered as the $q+q'-1$ first iterates of  a vertical line under a
rigid rotation of 
angle $\alpha\in ]\frac{p}{q},\frac{p'}{q'}[$. 
\end{prop}

In other words, to prove the  line translation theorem, it is enough to
find an essential  {topological line} $\gamma$ in $\AA$, and a lift of
$\gamma$ 
which   is   disjoint   from   its   images   under   $\wt\phi$   and
$\wt\psi$.  The  analogue  of proposition~\ref{p.arithmetics}  in  the
context of homeomorphisms  of the circle is well-known. The
proof  of  the proposition  relies  on arithmetical  properties of  Farey
intervals.  The  reader  can  find   a   proof 
in~\cite[appendix A]{beguin-crovisier-leroux-patou} (the proof is written in 
the context of the closed annulus, but also works 
in the infinite annulus setting).

\subsection{Proof of the line translation theorem} 
\label{ss.proof-arc-translation}

The closures of the rotation sets of the homeomorphisms $\wt
\phi=T^{-p}\circ\wt 
h^q$ and $\wt\psi=T^{p'}\circ\wt h^{-q'}$ are included respectively in
$]0,\frac{1}{q'}[$ and $]0,\frac{1}{q}[$ (see
proposition~\ref{p.usual-properties}). 
  In particular, according to
proposition~\ref{p.conjugue-a-translation},      the     homeomorphism
$\wt\phi$ is conjugate to a 
translation, and  thus, the quotient $\AA':=\RR^2/\wt\phi$  is an open
annulus. We  denote by $\pi'$ the natural projection of  $\RR^2$ onto
$\AA'$. 

Since   $\wt\phi$  and   $\wt\psi$  commute,   $\wt\psi$   induces  a
homeomorphism $\psi'$ of the open annulus $\AA'$.  The lift $\wt\psi$
of  $\psi'$ is  fixed  point  free. The next task is to check that
$\psi'$ satisfies the second hypothesis of theorem~(\ref{t.Guillou-Sauzet}). 

\bigskip

\noindent\emph{Claim 1. No point of the annulus $\AA'$ is wandering under the iteration of $\psi'$.} 

\begin{proof}
We  shall prove  that a dense set of points of $\AA'$  are
recurrent for the homeomorphism $\psi'$; the claim will follow. 

Poincar\'e recurrence theorem implies that a dense set of points of
$\AA$ are recurrent for $h$; this set lifts to a dense set in $\RR^2$,
which again projects to a dense set in $\AA'$. We prove that this last
set consists of recurrent points for $\psi'$. 

Consider a point  $\wt x\in\RR^2$, such that the  point $x:=\pi(\wt x)$ in $\AA$ 
is   recurrent   for  $h$.  
   There exists  two  sequences   of  integers
$(i_n)_{n\in\NN}$  and $(j_n)_{n\in\NN}$, such  that 
$j_{n}\rightarrow +\infty$ and
$T^{-i_n}\circ\wt
h^{j_n}(\wt x)\rightarrow\wt  x$
when $n$ goes  to $+\infty$. For every $n$,  we set 
$$k_n:=j_np'-i_nq'  \ \ \mbox{ and } \ \        l_n:=j_np-i_nq,$$
       so        that       
$$T^{-i_n}\circ \wt h^{j_n}=\wt\psi^{l_n}\circ\wt\phi^{k_n}.$$
  Hence,  
 $\wt\psi^{l_n}\circ\wt\phi^{k_n}(\wt x)\rightarrow\wt x$ when $n$
goes   to  $+\infty$,
which implies that 
${\psi'}^{l_n}(x)\rightarrow x$.
Moreover, since $\t \phi$ is conjugate to a translation, it has no
  recurrent point, so $l_{n}$ cannot be equal to zero for $n$ large
  enough 
(since $\wt\psi^{l_n}\circ\wt\phi^{k_n}(\wt x)\rightarrow\wt x$, this
  would imply $k_{n}=0$ for large $n$; since we also have $j_n=k_nq-l_nq'$,
and $j_n\rightarrow +\infty$, this would be a contradiction.). 
  Thus the point  $x':=\pi'(\wt  x)$  is
recurrent for the homeomorphism $\psi'$.  This complete the proof of claim~1.
\end{proof}

We are now in a position to  apply  Guillou-Sauzet-Le Calvez
theorem~\ref{t.Guillou-Sauzet}; 
it  provides  us  with a  Brouwer  line 
$\Gamma$ for $\wt\psi$, such that the projection $\gamma'$ of $\Gamma$
in the annulus $\AA'= \RR^2 / \t \phi$ is an essential
{topological line}. This implies that $\Gamma$ is also a Brouwer line
for $\t \phi$. 
According  to proposition~\ref{p.arithmetics},  we are  left  to prove
that the projection $\gamma$ of  the line $\Gamma$
\emph{in the original annulus $\AA=\RR^2/T$} is again an essential
{topological line}. 

\bigskip

\noindent\emph{Claim    2.    The    lines    $\wt\psi(\Gamma)$    and
  $\wt\phi(\Gamma)$ belongs to the same connected component of
  $\RR^2\setminus\Gamma$.}  

\begin{proof}
We choose
an orientation of $\Gamma$ in such a way that $\wt\phi(\Gamma)$ is on
the  right of  $\Gamma$ (see
subsection~\ref{ss.Brouwer-lines}).  For every  $k,l\in\ZZ$,
the     line
$\wt\phi^k\circ\wt\psi^l(\Gamma)$ is endowed with  
the image by $\wt\phi^k\circ\wt\psi^l$
of the orientation  of
$\Gamma$. We  denote by $U$ be  the connected open region of
$\RR^2$ bounded by the lines $\Gamma$ and $\wt\phi(\Gamma)$. 

We argue by contradiction: we  assume that $\wt\psi(\Gamma)$ is on the
left of  $\Gamma$, or equivalently, that  $\wt\psi^{-1}(\Gamma)$ is on
the  right  of $\Gamma$.  Under  this  assumption, the  homeomorphisms
$\wt\phi$  and $\wt\psi^{-1}$  are  both ``pushing  the line  $\Gamma$
towards the right''.  Hence, for every pair of positive integer $(k,l)$,
the  region   $\wt\phi^k\circ\wt\psi^{-l}(U)$  is  on   the  right  of
$\wt \phi(\Gamma)$, and thus is disjoint from $U$. 

According to Le Calvez theorem~\ref{th.mean-rotation-number},
almost every  point of the annulus $\AA$  is recurrent under $h$ 
and has a well-defined rotation number. Thus we can find a point $\wt x$ in
$U$ and some positive integers $m,n$ such that the point $\wt h^m\circ
T^{-n}(\wt x)$ is in $U$  and such that $n/m$ belongs to
$]p/q,p'/q'[$.  We have
$$
\wt h^m\circ T^{-n}=\wt\phi^{k}\circ\wt\psi^{-l},
$$ 
$$
\mbox{with } k=mp'-nq' \mbox{ and } l=-mp+nq.
$$
  Since $n/m$ is  in the Farey  interval $]p/q,p'/q'[$,
the  integers  $k=mp'-nq'$ and  $l=-mp+nq$  are  positive. Hence,  the
region $\wt h^m\circ T^{-n}(U)$ is disjoint from the region $U$. But
 this is absurd,
since the  point $h^m\circ  T^{-n}(\wt x)$ is  in the  intersection of
these two regions.
\end{proof}

\noindent\emph{Claim 3.  The line $\Gamma$ is a  Brouwer line for $T$.
Furthermore, let $V$   be the
  connected open region  of $\RR^2$  bounded by  the  lines $\Gamma$  and
  $T(\Gamma)$, and $\adhe(V)=\Gamma \cup V \cup T(\Gamma)$. Then $\adhe(V)$  
  is    a   fundamental   domain    for   the   covering map 
  $\RR^2\rightarrow\AA=\RR^2/T$.}  

\begin{proof}
By claim 2, both homeomorphisms $\wt\phi$ and $\wt\psi$ ``push the
line   $\Gamma$   towards   right''.   Hence,   given   four   integers
$k,l,k',l'\in\ZZ$,   such   that   $k<k'$   and   $l<l'$,   the   line
$\wt\phi^{k'}\circ\wt\psi^{l'}(\Gamma)$  is strictly  on the  right of
the line $\wt\phi^k\circ\wt\psi^l(\Gamma)$ (we call this ``property $(\star)$'').

In   particular,   $T(\Gamma)=\wt\phi^q\circ\wt\psi^{q'}(\Gamma)$   is
strictly on the right of $\Gamma$, and $T^{-1}(\Gamma)$ is strictly on
the left of $\Gamma$. Therefore, $\Gamma$ is a Brouwer line for $T$. 

We are  left to prove  that the iterates  of $\adhe(V)$  under $T$
cover   the   whole   plane,  \emph{i.e.}    that   $\bigcup_{k\in\ZZ}
T^k(\adhe(V))=\RR^2$.  As above,  we  denote  by $U$  the  connected
open region of  $\RR^2$ bounded by  the lines $\Gamma$ and
$\wt\phi(\Gamma)$. Since the 
projection  of  $\Gamma$ in  the  annulus  $\AA'=\RR^2/\wt\phi$ is  an
essential simple  {line}, $\adhe(U)=\Gamma \cup U \cup \wt\phi(\Gamma)$
  is a fundamental  domain for  the covering map
$\RR^2\rightarrow\AA'$,        and         thus, we have
$$
\bigcup_{k\in\ZZ}\wt\phi^k(\adhe(U))=\RR^2.
$$
 According to property $(\star)$,
for every $n>0$, the line
$T^{-n}(\Gamma)=\wt\phi^{-nq}\circ\wt\psi^{-nq'}(\Gamma)$  is  on  the
left    of   the   line    $\wt\phi^{-n}(\Gamma)$, and    the   line
$T^{n}(\Gamma)=\wt\phi^{nq}\circ\wt\psi^{nq'}(\Gamma)$ is on the right
of  the   line  $\wt\phi^{n}(\Gamma)$ (remember that $q$ and $q'$ are
 greater than $1$). 
   Now  observe  that   the  set
$\bigcup_{k=-n}^{n-1} T(\adhe(V))$  is the region situated between the lines
$T^{-n}(\Gamma)$      and     $T^n(\Gamma)$,      and      the     set
$\bigcup_{k=-n}^{n-1}\wt\phi(\adhe(U))$  is the  region situated
 between the
lines    $\wt\phi^{-n}(\Gamma)$   and   $\wt\phi^n(\Gamma)$.    As   a
consequence, for every $n>0$, we have 
$$
\bigcup_{k=-n}^{n-1} T(\adhe(V))\supset \bigcup_{k\in\ZZ=-n}^{n-1}\wt\phi^k(\adhe(U)),
$$
                and               thus
$\bigcup_{k\in\ZZ}T^k(\adhe(V))=\RR^2$.
  This  completes  the  proof  of  the claim.  
\end{proof}

\noindent\emph{Claim 4.  The line $\Gamma$ projects in $\AA$ to an essential 
 {line} $\gamma$}

\begin{proof}
What remains to be proved is that, with respect to  the translation
$T: (x,y) \mapsto (x+1,y)$, the Brouwer line $\Gamma$ 
 is equivalent to the ``trivial'' Brouwer line 
$\Gamma_0:=\{0\} \times \RR$. That is, that $\Gamma$ is proper in $\AA$.
For that, it suffices to construct a homeomorphism
$G$ of the plane that commutes with $T$, and such that $G(\Gamma_0)=\Gamma$.
This is very classical, as we have already mentioned in 
paragraph~\ref{ss.Brouwer-lines}. By Schoenflies theorem, there exists
a homeomorphism $G$ from
the band $[0,1] \times \RR$ onto the region $\adhe(V)$, such that
$$
T \circ G_{|\{0\} \times \RR} =  G_{|\{1\} \times \RR} \circ T.
$$
 Then we extend $G$ by conjugacy, 
that is, we set
$$
G(p+\alpha,t)=T^p(G(\alpha,t))
$$
for any real number $t$, any integer $p$ and any number $\alpha$ 
between $0$ and $1$.
The map $G$ is continuous. It is one-to-one (because $\Gamma$ and $V$ are disjoint 
from their iterates under $T$). It is onto (because of claim 3).
Clearly, $G$ is an open map; hence, it is an homeomorphism.
\end{proof}

This completes the proof of the line translation theorem.

\section{Closure of the conjugacy class of a pseudo-rotation}
\label{s.closure}

Recall that theorem~\ref{t.closure-conjugacy-class} states that, for
any pseudo-rotation $h:\AA\rightarrow\AA$ of irrational angle $\rho$,
the rigid rotation of $(x,y)\mapsto (x+\rho,y)$ is in the closure (for
the compact open topology) of the conjugacy class of $h$.
 A similar
result was proved in~\cite[corollary
0.2]{beguin-crovisier-leroux-patou} in the compact annulus 
setting. Actually, the proof given
in~\cite[section 5]{beguin-crovisier-leroux-patou} applies to the
open annulus  
setting, with the following modifications:  
\begin{itemize}
\item[--] replace the notion of \emph{essential simple arc} used
in~\cite{beguin-crovisier-leroux-patou} by the notion of
\emph{essential topological line} defined in the present article,
\item[--] instead of using the \emph{arc translation theorem}
of~\cite{beguin-crovisier-leroux-patou}, use the \emph{line translation
theorem} of the present article.
\end{itemize}

\section{Integrability of the displacement function}
\label{s.integrabilite}

The aim of this section is to show that, under the hypothesis of Le Calvez theorem~\ref{th.mean-rotation-number}, up to a suitable change of coordinates, the horizontal displacement function is bounded, and hence integrable (the horizontal displacement function has been defined in paragraph~\ref{ss.morphisms}).

\subsection{Statements}
\label{ss.statements}

\begin{prop}[Integrability of the displacement function]
\label{p.integration}
Consider a homeomorphism $h \in \homeo^+_{\leb}(\AA)$.  Assume that the
set $\rot_{\fix}(\t h)$  of rotation numbers of the fixed points of $h$ is bounded (for some lift $\t h$). 

Then  there exists  $g \in  \homeo^+_{\leb}(\AA)$,  such that
the horizontal displacement  function  $r$  of  any  lift  $\t h_{1}$ of  the  homeomorphism
$h_1=g\circ h \circ g^{-1}$ is bounded.
\end{prop}

Note that as a consequence of Birkhoff ergodic theorem, the mean rotation number of $\t h_{1}$  is equal to the integral of $r$ over the annulus $\AA$.
As a classical consequence, we get a more geometrical definition.

\begin{prop}
\label{p.geometric-interpretation}
Let $h_{1} \in \homeo_{\leb}^+(\AA)$, and $\t h_{1} : \t \AA \rightarrow \t \AA$ be a lift of $h_{1}$. Suppose that the horizontal displacement function $r$ of $\t h_{1}$ is bounded.
Then the mean rotation
number of $\t h_{1}$ is equal to the algebraic area (for the
lift of the Lebesgue probability measure on $\AA$) of the region of
$\t \AA=\RR\times\RR$   situated  between  any   vertical  line
$\t D=\{\theta\}\times\RR$ and its image $\t h_1(\t D)$.  
\end{prop}


In  view to proposition~\ref{p.integration}, it  seems natural  to hope
that (under  suitable assumptions) the mean rotation  number ``defines a
morphism'', as in the case of the compact annulus (see~\ref{ss.morphisms}).
 For example, the following question may be asked. 

\begin{ques}\label{q.morphism}
Let $f$, $g$  be two homeomorphisms of the  annulus, which are isotopic
to the identity and preserve the Lebesgue probability measure. Consider some lifts
$\wt f,\wt g$  of $f,g$, and assume that the  mean rotation numbers of
$\wt f$, $\wt g$ and $\wt f\circ\wt g$ are well-defined.  

Is the mean  rotation number of $\wt f\circ\wt g$ equal  to the sum of
the mean rotation numbers of $\wt h$ and $\wt g$~? 
\end{ques}

We briefly explain the idea of the proof of
proposition~\ref{p.integration}. The easy case is when the closure of
the rotation set of $\t h$ is contained in some interval $]p, p+1[$
with $p\in\ZZ$ (e.g. when $h$ is an irrational pseudo-rotation).
In this case, since $]p,p+1[$ is a Farey interval, we can directly apply the line translation theorem~\ref{th.arc-translation}, and we get an essential topological line in $\AA$ which is disjoint from its image under $h$. The conservative version of Schoenflies theorem gives a $g \in \homeo^+_{\leb}(\AA)$ that maps this topological line on the straight line $\{0\} \times \RR$. The conjugated homeomorphism $g h g^{-1}$ now maps this straight line off itself, and we see easily that the horizontal displacement function of any lift is bounded. In the general case, we will use this easy case by considering intermediate coverings.

\subsection{Rotation numbers for intermediate coverings}
\label{ss.intermediate}
As usual, take $h \in \homeo^+_{\leb}(\AA)$ and $\t h: \t \AA \rightarrow \t \AA$ a lift of $h$.
Remember that $T$ denotes the covering translation of $\t \AA$ (which commutes with $\t h$).
Given an integer $q \geq 2$, we may consider the intermediate covering
$\AA_q=\t \AA/T^q$, which is again an annulus.  
The homeomorphism $\t h$ induces a homeomorphism $h'$ of $\AA_{q}$.
In addition to the previously defined notions of rotation numbers of $\t h$ as a lift of $h$, one can consider the rotation numbers of $\t h=\t h'$ \emph{as a lift of $h'$}.
These numbers are linked in the following way.
If $z$ is a recurrent point for $h$, and $z'$ is any lift of $z$ in $\AA_{q}$, then one easily proves that $z'$ is a recurrent point for $h'$. 
If $z$ has a well-defined rotation number $\rho(z, \t h)$ under $h$, then the rotation number of $z'$ under $h'$ is also well defined and equal to $\frac{1}{q} \rho(z, \t h)$.

\subsection{Proofs}
The core of the proof of proposition~\ref{p.integration} is contained
in the proposition given below. 
We use the notations of the previous paragraph.
Assume that the closure of the rotation set of  $\t h'=\t h$ as a lift
of $h'$ is 
contained in the Farey interval $]0,1[$. Then we can apply the line
translation theorem~\ref{th.arc-translation}, which provides us with
an essential topological line $\gamma'$ of $\AA_{q}$, which is
disjoint from its image $h'(\gamma')$. Note that in general, the
projection of $\gamma'$ in $\AA$ is not a topological line (it may
have self-intersections). 

\begin{prop}
\label{p.commute}
We can choose the topological line $\gamma'$ so that its projection in
$\AA$ is again a  topological line. 
\end{prop}

We will also call {\em essential topological  line in $\RR^2$} an
oriented simple curve $\Gamma: \RR \rightarrow \RR^2$ such that the
second coordinate of $\Gamma(t)$ tends to $-\infty$ (resp. $+\infty$)
when $t$ tends to $-\infty$ (resp. $+\infty$).
 Remember that we denote by $R(\Gamma)$ (resp. $L(\Gamma)$) the connected component of
 $\RR^2\setminus \Gamma$  on the right (resp. on the left) of $\Gamma$. 
If $\Gamma_1$ and $\Gamma_2$ are two essential topological lines in
$\RR^2$, we write $\Gamma_1 \leq \Gamma_2$  when $\Gamma_2$ is contained in
$R(\Gamma_1)$; we write  $\Gamma_1 < \Gamma_2$ if $\Gamma_1 \leq \Gamma_2$ and the two lines are disjoint. 

\begin{lemmanota}
\label{l.vee}
Let $\Gamma_1$ and $\Gamma_2$ be two essential topological lines in 
$\RR^2$, and let $U$ be the unique connected component of the set 
$L(\Gamma_1)\cap
L(\Gamma_2)$ which contains half lines of the form
$]-\infty,a[\times \{b\}$. Then the boundary
of $U$ is an essential topological
line in $\RR^2$, that we denote by $\Gamma_1\vee\Gamma_2$.
\end{lemmanota}

The     proof  of lemma~\ref{l.vee}  is     similar     to    that     of    lemma     3.2
in~\cite{beguin-crovisier-leroux-patou} and uses a classical result by
B. Ker\'ekj\'art\'o (\cite{kerekjarto2}).  

\begin{rema}
\label{r.vee}
Let $\Gamma_1$, $\Gamma_2$, $\Gamma_3$  be three essential topological lines
in $\RR^2$. The following
properties are immediate consequences of the definition of the line
$\Gamma_1\vee\Gamma_2$.

\begin{itemize}
\item[\textbf{\emph{(i)}}] The line $\Gamma_1\vee\Gamma_2$ is 
included in the union of the lines $\Gamma_1$ and $\Gamma_2$. Hence, 
if      $\Gamma_3<\Gamma_1$     and      $\Gamma_3<\Gamma_2$,     then
$\Gamma_3<\Gamma_1\vee\Gamma_2$. 
\item[\textbf{\emph{(ii)}}] The sets $R(\Gamma_1)$ and 
$R(\Gamma_2)$ are  included in the  set $R(\Gamma_1\vee\Gamma_2)$.  In
other   words,    we   have   $\Gamma_1\vee\Gamma_2\leq\Gamma_1$   and
$\Gamma_1\vee\Gamma_2\leq \Gamma_2$. 
\end{itemize}
\end{rema}

\begin{proof}[Proof of proposition~\ref{p.commute}]
By theorem~\ref{th.arc-translation}, there  exists an essential topological
line $\gamma_0$ of $\AA_q$ which is disjoint from its image $h'(\gamma_{0})$. We
consider  some  lift  $\Gamma_0$   of  $\gamma_0$  to  $\RR^2$.  Since
$\gamma_0$ is simple  in $\AA_q$, the arc $\Gamma_0$  is disjoint from
$T^q(\Gamma_0)$. Note that  since the rotation set of  $\t h'=\t h$ as a
lift  of $h'$ is  contained in  $]0,1[$, we  have $T^{-q}(\Gamma_0)<\t
h^{-1}(\Gamma_0)<\Gamma_0$. 

Now, we  choose some essential topological lines $\Gamma_1,\dots,\Gamma_{q-1}$ in
$\RR^2$ such that 
$$T^{-q}(\Gamma_0)<\Gamma_{q-1}<\Gamma_{q-2}<\dots<\Gamma_{1}<\Gamma_0.$$ 
Consider the essential topological line 
$$\Gamma=\Gamma_0\vee        T(\Gamma_{1})\vee      \dots\vee
T^{q-1}(\Gamma_{q-1})=\bigvee_{i=0}^{q-1}T^i(\Gamma_i).$$
For         every        $i\in\{0,\dots,q-2\}$,       we have
$T^{i+1}(\Gamma_{i+1})<T^{i+1}(\Gamma_i)$  (by  definition
of the $\Gamma_i$'s) and
$\Gamma\leq T^{i+1}(\Gamma_{i+1})$  (by definition of  $\Gamma$ and by
item~(ii) of 
remark~\ref{r.vee}). Hence for every $i\in\{0,\dots,q-2\}$, we get 
$$\Gamma<T^{i+1}(\Gamma_i).$$ 
Moreover,    we    have    $\Gamma_0<T^q(\Gamma_{q-1})$  and
$\Gamma\leq\Gamma_0$.  Hence
$$\Gamma<T^q(\Gamma_{q-1}).$$ 
Finally, using item~(i) of remark~\ref{r.vee}, we get 
$$\Gamma<\bigvee_{i=0}^{q-1}T^{i+1}(\Gamma_i)=T(\Gamma).$$
In particular, $\Gamma$ is disjoint from its image under $T$. 
Moreover,      we  may     assume    that      the      lines
$\Gamma_{1},\dots,\Gamma_{q-1}$ were chosen such that 
$$\t h^{-1}(\Gamma_0)< \Gamma_{1},\dots,\Gamma_{q-1}<\Gamma_0.$$
Using the  definition of $\Gamma$ and  remark~\ref{r.vee}, this easily
implies that
$$\Gamma<\t h'(\Gamma).$$ 
Similarly, since $\t h \circ T^{-q} (\Gamma_{0}) < \Gamma_{0}$,     we   may      assume      that      the      lines
$\Gamma_{1},\dots,\Gamma_{q-1}$ were chosen such that 
$$\t       h\circ       T^{-q}(\Gamma_0)<\Gamma_{1},\dots,\Gamma_{q-1}
<\Gamma_0.$$
This easily implies that
$$\t h'(\Gamma)<T^{q}(\Gamma).$$ 

Let $\gamma'$ be the projection of $\Gamma$ in the annulus $\AA_{q}$. Since $\Gamma < \t h(\Gamma) < T^{q}(\Gamma)$, the curve $\gamma'$ is an essential topological line in $\AA_{q}$ which is disjoint from its image $h'(\gamma')$. Furthermore, since $\Gamma < T(\Gamma)$, the projection of $\gamma'$ in the annulus $\AA$ is again simple, thus it is an essential topological line.
\end{proof}

We are now able to prove proposition~\ref{p.integration}.

\begin{proof}[Proof of proposition~\ref{p.integration}]
Let $h\in\homeo^+_\leb(\AA)$, and $\t h$  be a lift of $h$.
Assume that the set of rotation numbers of the fixed points of $h$ is bounded.
Applying theorem~\ref{t.PBFLC} to the homeomorphisms $T^{-p} \circ \t h$, we see that the rotation set of $\t h$ is also bounded. Up to a change of lift, we may assume that $\rot(\t h)$ is included in an interval
$[1,q-1]$ for some integer $q$.

 Consider the homeomorphism $h'$ induced by $\t h$ on the intermediate
covering $\AA_q$.
The rotation set of
$\t  h$, seen as a lift of $h'$, is included in $[\frac{1}{q},\dots,\frac{q-1}{q}]$ (see paragraph~\ref{ss.intermediate}).
Hence, by proposition~\ref{p.commute}, there exists 
an essential topological line $\gamma$ in $\AA$  and a lift $\gamma'$ of $\gamma$ in $\AA_q$
which is disjoint from its image $h'(\gamma')$.

Using the conservative Schoenflies theorem, we get some homeomorphism $g\in\homeo^+_\leb(\AA)$
which sends $\gamma$ on the vertical line $D=\{0\}\times \RR$.
Take any lift $\t  g$ of $g$ to $\RR^2$ and denote by $g'$ the induced map on $\AA_q$.
We consider the conjugates $h_1=g\circ h\circ g^{-1}$, $h'_1=g'\circ h'\circ (g')^{-1}$
and $\t  h_1=\t  g\circ \t  h\circ \t  g^{-1}$. 
Consider some lift $\t D$ of $D$ in $\RR^2$, and its projection $D'$ in $\AA_{q}$.
Since $h'(\gamma')$ is disjoint from $\gamma'$, the line $D'$ is disjoint from $h'_{1}(D')$.
Thus the image $\t h_{1}(\t D)$ is between $T^{qk}(\t D)$ and $T^{qk+q}(\t D)$ for some integer $k$.
This easily implies that the image of the displacement function $r=p_1\circ \t  h_1-p_1$
is bounded.
The proof is complete.
\end{proof}

\begin{proof}[Proof of proposition~\ref{p.geometric-interpretation}]
We may assume that $\t D$ is the vertical line $\{0\}\times \RR$. We also note that
if the proposition is satisfied for some lift $\t h_1$ of $h_1$, then it is true for any other lift.
Since the displacement functions of the lifts of $h_1$ are bounded, one may choose a lift $\t h_1$
and an integer $q\geq 1$ such that $\t D< \t h_1(\t D)< T^q(\t D)$.
Hence, there exists a vertical line $\t D_1$ between $\t D$ and $T^q(\t D)$ such that the area
of the strip bounded by $\t D$ and $\t h_{1}(\t D)$ is equal to the area of the strip
between $\t D$ and $\t D_1$.
We will work on the intermediate covering $\AA_q$ with the homeomorphism $h_{1}'$ induced by $\t h_{1}$.
In $\AA_q$, the projection of $\t D_1$ is a vertical line $D'_1$ of the form
$\{\theta\}\times \RR$ and the projection of $\t D$ is a vertical line $D'$.

By the conservative Schoenflies theorem, there exists
some homeomorphism $\psi\in\homeo^+_\leb(\AA_q)$ which is the identity on $D'$ and which maps
the line $h'_{1}(D')$ on $D'_1$: one can require moreover that each point $(0,t)$ in $D'$
is mapped by $\psi \circ h'_1$ on the point $(\theta,t)$ in $D'_1$. We denote by $\t \psi$ the lift of $\psi$
to $\t \AA_q=\t \AA$ which fixes $\t D$. Note that the displacement function of $\t \psi$ is bounded.
We also introduce the translation $R_{\theta}\colon \t \AA\to \AA$
by $\theta$. It is the lift of the rotation $R'$ of $\AA_q$ with angle $\frac{\theta}{q}$.

Since $\t \psi$ and $R_{\theta}^{-1}\circ \t\psi \circ \t h_1$ are the identity on $\t D$, 
their rotation sets (as lift of homeomorphism on $\AA_q$) are $\{0\}$ and
their mean rotation number are $0$.
By the morphism property (proposition~\ref{p.morphism-property}),
one deduces that the mean rotation number of $\t h_1$ (as lift of homeomorphisms of $\AA_q$)
is equal to $\frac{\theta}{q}$. As it is explained at paragraph~\ref{ss.intermediate}, this implies that the mean rotation
number of $\t h_1$, as lift of the homeomorphism $h_1$ of $\AA$, is equal to $\theta$.
By construction, $\theta$ is also the area between $\t D$ and its image $\t h_1(\t D)$
in $\t \AA$. This concludes the proof.
\end{proof}


\section{Periodic orbits in one-parameter families}
\label{s.familles}
In this section, we consider perturbations of conservative homeomorphisms by  hamiltonian vector flows, and we prove a  general statement that implies theorem~\ref{t.perturbation-pseudo-rotation}.

\subsection{General statement and some consequences}
 Let $X$ be a smooth vector field on $\AA$ such that
\begin{enumerate}
\item $X$ is bounded (in the coordinates system $\AA = \SS^1 \times \RR$);
\item the flow $(\Phi^t)_{t \in \RR}$ generated by $X$ preserves the Lebesgue probability measure on $\AA$.
\end{enumerate}

The vector field $X$ lifts to a vector field $\tilde X$ on $\tilde \AA$, which in turn generates a flow $(\tilde \Phi^t)_{t\in \RR}$ which is a lift of the flow $(\Phi^t)_{t\in\RR}$.
The rotation number of the Lebesgue probability measure 
for $\tilde \Phi^t$ is well defined. Since the displacement function is bounded, the morphism
property of proposition \ref{p.morphism-property}
is satisfied and the map $t\mapsto \rho(\leb,\tilde \Phi^t)$ is continuous.
Hence, we have 
$$
\rho(\leb, \tilde \Phi^t)= t . \rho(\leb,\tilde \Phi^1).
$$
Our favourite example is of course the family of Euclidean rotations on $\SS^2$ (given on $\AA=\SS^1 \times \RR$ by the constant vector field $X=(1,0)$).
The following statement implies at once theorem~\ref{t.perturbation-pseudo-rotation}.
\begin{theo}
\label{t.hamiltonian-perturbation}
Let $X$ be a vector field as above, and suppose 
$$
\rho(\leb,\tilde \Phi^1) \neq 0.
$$
Let $h \in \homeo_\leb^+(\AA)$.
Then there  exists
arbitrarily  small values  of  $t$  such that  $h\circ  \Phi^{t}$ has  a
periodic orbit.
\end{theo}

Let us first state and prove two interesting corollaries of this theorem.
\begin{coro}
\label{c.one-parameter-families}
Given any homeomorphism $h\in\homeo^+_\leb(\AA)$, the set 
$$D(h):=\{t\in\RR| h \circ \Phi^t \mbox{ has a periodic orbit}\}$$
is dense in $\RR$.
\end{coro}
For each $r\in [0,\infty]$, we consider the  space $\mathrm{Diff}^{r,+}_\mathrm{Leb}(\AA)$
of $C^r$  diffeomorphisms
 of the annulus $\AA$ that preserves
the  Lebesgue probability measure, endowed either  with the  compact-open  or the
Whitney topology.
\begin{coro}
\label{c.meager}
The space $\mathrm{Diff}^{r,+}_\mathrm{Leb}(\AA)$ contains an open and dense subset of diffeomorphisms having periodic orbits.
In particular,  the set of irrational pseudo-rotations is meagre.
\end{coro}

\begin{proof}[Proof     of    corollary~\ref{c.one-parameter-families}
  assuming theorem~\ref{t.hamiltonian-perturbation}] 
Fix  $t_0\in\RR$.  Applying  theorem~\ref{t.hamiltonian-perturbation} to  the
homeomorphism $h\circ \Phi^{t_0}$,  we get arbitrarily small numbers
$\eta$  such that  the homeomorphism  $h\circ \Phi^{t_0+\eta}=h\circ
\Phi^{t_0}\circ \Phi^\eta$ has a periodic orbit.  
\end{proof}

\begin{proof}[Proof      of      corollary~\ref{c.meager}     assuming
  theorem~\ref{t.hamiltonian-perturbation}] 
Consider    the     set    $\cU$     of    all    the     elements    of
$\mathrm{Diff}^{r,+}_\mathrm{Leb}(\AA)$   that  have   at   least  one
hyperbolic  periodic  orbit.  Note  that  $\cU$ is  an  open  subset  of
$\mathrm{Diff}^{r,+}_\mathrm{Leb}(\AA)$ (since  hyperbolic  periodic
orbits   are   persistent) for the compact-open toplogy, and so for the Whitney topology.
Consider  an   element   $h$   of
$\mathrm{Diff}^{r,+}_\mathrm{Leb}(\AA)$, and an open  neighbourhood $\cV$ of $h$ in $\mathrm{Diff}^{r,+}_\mathrm{Leb}(\AA)$, for the Whitney topology.
There exists  $h'\in\cV$ that has  a periodic orbit: 
this follows  from theorem~\ref{t.hamiltonian-perturbation} since there exists a smooth vector field
$X$ of $\AA$ with compact support such that the time-one map $\Phi^1$ has a non-zero rotation number. We can perturb
$h'$ in order 
to get  a diffeomorphism $h''\in\cV$ that  has  a
hyperbolic  periodic  orbit.   This   proves  that  $\cU$  is  dense  in
$\mathrm{Diff}^{r,+}_\mathrm{Leb}(\AA)$ for the Whitney topology, hence also for the compact-open topology.
\end{proof}

\subsection{Idea of the proof}
Let   us  explain   briefly   the   idea   of  the   proof   of
theorem~\ref{t.hamiltonian-perturbation}, as developed in paragraph~\ref{ss.general}. We suppose $h$ is an irrational pseudo-rotation, with
 angle $\alpha$ (otherwise there is  nothing to prove).
 There  are two disjoint
cases.  Either  the  rotation  sets  of  the homeomorphisms  $h  \circ \Phi^t$ ``explode'' (that is, there exists arbitrarily 
small  values  of $t$  for  which  the  rotation set  of  $h\circ
\Phi^t$  contains numbers  arbitrarily  far from  $\alpha$), or  the
rotation set of  $h \circ \Phi^t$ is uniformly  bounded for $t$
close to $0$. Surprisingly, the first case is the easiest: because
of  the  lower  semi-continuity  property  of the  rotation  set (see below),  the
rotation  set   of  $h\circ  \Phi^t$  must  also contain   some  numbers
arbitrarily    close     to    $\alpha$;    so     we    can    apply
Poincar\'e-Birkhoff-Franks's  theorem  to  get  a  periodic  orbit  of
rotation number close to $\alpha$.  In the second case, we use 
proposition~\ref{p.integration} that allows us to suppose that  the horizontal
displacement  of  $h$ is  bounded.  Approximating  the maps $\Phi^t$  by
compactly supported maps and using a continuity property (see below),
 we see that the ``morphism property'' holds:
$\rho(\leb, \tilde h \circ \tilde \Phi^t)=\rho(\leb,\tilde h)+\rho(\leb,\tilde \Phi^t)$. When
$\alpha+t$ is rational, this again gives rise to periodic orbits.

\subsection{The compactly supported case}
Before addressing the general issue, it is useful to deal with a restricted problem.
\emph{We first prove theorem~\ref{t.hamiltonian-perturbation} assuming that the vector field $X$ is compactly supported in $\AA$.} The following argument is essentially due to J. Franks ( \cite{franks3}).
\begin{proof}[Proof (compactly supported case)]
If $h$ is not an irrational pseudo-rotation, then it has a periodic orbit (proposition~\ref{p.no-periodic-orbit}), hence we can take $\theta=0$, and there is nothing to prove. So from now on we assume that $h$ is an irrational pseudo-rotation.

We fix a lift $\t h$ of $h$. We denote by $\alpha$ the rotation number
of  $\t h$. 
Since $h$ is a pseudo-rotation, its rotation set is certainly bounded. Thus we can apply proposition~\ref{p.integration}\footnote{Note that this is the easy case of the proposition, as explained at the end of paragraph~\ref{ss.statements}.}: by performing a change of coordinates given by a homeomorphism $g$, we may assume that the horizontal displacement function of $\t h$ is bounded. Note that the change of coordinates does not affect the fact that the flow $(\Phi^t)_{t \in \RR}$ is compactly supported\footnote{The conjugated flow is not smooth anymore, which will not do any harm.}.

We now deal with $\t h$ and $\t \Phi^t$ having bounded (integrable) horizontal displacement functions. Thus we have the morphism property (proposition~\ref{p.morphism-property}):
$$
\rho(\leb,\t   h    \circ   \t \Phi^t)=\rho(\leb ,\t   h    )+   \rho(\leb,
\t \Phi^t)= \alpha+t . \rho(\leb, \t \Phi^1).
$$
By hypothesis, $\rho(\leb, \t \Phi^1)$ is non null, so there exists arbitrarily small values of $t$
such that the mean rotation number $\rho(\leb,\t   h    \circ   \t \Phi^t)$ is rational.
For any such value, $h \circ \Phi^t$ is not an irrational pseudo-rotation, so it must have periodic orbits according to proposition~\ref{p.no-periodic-orbit}.
This solves the compactly supported case.
\end{proof}

\subsection{Some continuity results by P. Le Calvez}
\label{ss.continuity}
We need some more tools before coping with the general case.
Le  Calvez  has  proved  the  following continuity  property  for  the
rotation number of the Lebesgue probability measure (see~\cite[theorem~2]{lecalvez}).
Remember that $\rot_{\fix}(\t   h)$ denotes the set of rotation numbers of the fixed points of $\t h$ (see section~\ref{ss.definition}).
\begin{theo}[P.~Le~Calvez, continuity of the mean rotation number]
\label{th.continuity}
Consider   a   sequence   $(h_n)_{n\in\NN}$  of   $\homeo^+_\leb(\AA)$
converging towards some homeomorphism $h\in\homeo^+_\leb(\AA)$ for the 
compact-open topology. Consider also a sequence $(\t  h_n)$ of
lifts converging towards a lift $\t  h$ of $h$. 

If  the sets  $\rot_{\fix}(\t   h_n)$  are uniformly  bounded,  then the  set
$\rot_{\fix}(\t   h)$ is  bounded and  the sequence  of rotation numbers
$\rho(\leb,\t   h_n)$  converges  towards 
the  rotation number  $\rho(\leb,\t  h)$.
\end{theo}

Actually,  the proof  of  the above  theorem  uses another  continuity
property proved in the same paper (see~\cite[proposition~3]{lecalvez}).
We consider a sequence of lifts $(\t  h_n)$ converging towards a lift $\t h$ as in the previous statement.
\begin{theo}[P.~Le~Calvez, lower semi-continuity of the rotation set]
\label{th.semicontinuity}
If  the closures  of  the  rotation sets  $\rot(\t   h_n)$   converge  towards   some   interval $[a,b]\subset [-\infty,+\infty]$, then the closure of the rotation set
$\rot(\t h)$ is contained in $[a,b]$. 
\end{theo}
There are some easy remarks in this footnote\footnote{In particular,  the closure  of the conjugacy  class of  an irrational
pseudo-rotation  with angle  $\alpha$  contains only  pseudo-rotations
with angle  $\alpha$. However, the  closure of the conjugacy  class of
some homeomorphism  whose rotation set is  not reduced to  a point may
contain  some  homeomorphism  whose   rotation  set  is  smaller.  For
instance, since  $\AA$ is open, it  is easy to build  an example where
each homeomorphism $\t  h_n$ has  rotation set equal to $[0,1]$ and
$h$  has  a rotation  set  reduced to  $\{0\}$.  This  shows that  the
rotation set is not upper semi-continuous. 

Theorem~\ref{th.semicontinuity} is  false if one does  not assume that
the $h_n$'s  preserve the Lebesgue probability  measure.
Note that  in the context of  the compact
annulus,  the  rotation set  is  always  upper semi-continuous  (even without
assuming that the Lebesgue measure is preserved).}.

\subsection{The general case}
\label{ss.general}
\emph{We now cope with the general case, without assuming that $X$ is compactly supported.}
\begin{proof}[Proof (general case)]
As before, it suffices to consider an irrational pseudo-rotation.
We use the notations introduced for the compactly supported case.
 We consider two disjoint subcases.
\paragraph{First  subcase: the rotation  set of  $\t h  \circ \t \Phi^t$
is not uniformly bounded for $t$ close to $0$.} 
More precisely, there exists a  sequence of numbers
$t_n\rightarrow 0$ such that the set $\mathrm{Rot}(\t h \circ
\t \Phi^{t_n})$ contains  a number $\beta_n$  with $\beta_n\ra +\infty$
or $\beta_n\ra -\infty$.  

According   to  the   lower  semi-continuity   of  the   rotation  set
(theorem~\ref{th.semicontinuity}), there  must exist another sequence
$(\alpha_n) \ra \alpha$  such that $\alpha_n\in\mathrm{Rot}(\t h \circ
\t \Phi^{t_n})$.  For  each value  of  $n$,  choose  a rational  number
$\alpha'_n$ strictly  between $\alpha_n$ and $\beta_n$, in  such a way
that   the   sequence   $(\alpha'_n)$   tends  to   $\alpha$.    Then,
Franks'       version       of       Poincar\'e-Birkhoff       theorem
(first part of theorem~\ref{t.PBFLC}) provides for  each $n$ a periodic
orbit  for the  map $h  \circ \Phi^{t_n}$  with rotation  number
$\alpha'_n$, and we are done.

\paragraph{Second case: the rotation set  of $\t h \circ \t \Phi^{t}$ is
  uniformly bounded for $t$ close to $0$.} 
We now assume that there exists a bound $M>0$ and
an angle $t_0$ such that 
$$\mathrm{Rot}(\t h \circ \t \Phi^t) \subset [-M, M]$$ 
for all $t\in [-t_0,t_0]$. 
In particular, the rotation numbers of all the fixed points of $h\circ \Phi^t$
are     included     in     $[-M,M]$     for     every
$t\in[-t_0,t_0]$. 

Let $(\phi_{s})_{s \in [1,+\infty]}$ be a continuous family of smooth  functions from $\AA=\SS^1\times\RR$ to $[0,1]$ such that
\begin{enumerate}
\item $\phi_{s}$ is equal to $1$ on the compact annulus $\SS^1 \times [-s,s]$;
\item  $\phi_{s}$ is equal to $0$ outside the compact annulus $\SS^1 \times [-2s,2s]$.
\end{enumerate}
For each value of $s$ we consider the vector field $X_{s}= \phi_{s}.X$.
We denote by $(\Phi_{s}^t)_{t \in \RR}$ the corresponding flow. The lifted flow (generated by the pullback on $\t \AA$ of the vector field $X_{s}
$) is denoted by $(\t \Phi_{s}^t)_{t \in \RR}$. For each finite value of $s$, the flow 
$(\Phi^t_{s})_{t \in \RR}$ is compactly supported; whereas it coincides with the original flow  
$(\Phi^t)_{ t\in \RR}$ for $s=+\infty$.

Since we have slowed down the flow, we have the following easy but crucial property:
\emph{for any positive $t$, for any point $\t x \in \t \AA$ and any $s \geq 1$ there exists a time $t'$ between $0$ and $t$ such that $\t \Phi^{t'}(\t x) = \t \Phi^{t}_{s}(\t x)$.} In particular, for every $t\in[-t_0,t_0]$ and every $s \geq 1$, each rotation number of some fixed point of $h \circ \Phi_{s}^t$
is equal to a rotation number of some fixed point of $h \circ \Phi^{t'}$ for $t'\in[-t_0,t_0]$. So the sets $\rot_{\fix}(\t h \circ \t \Phi_{s}^t)$ of rotation numbers of the fixed points are all included in $[-M,M]$.

With this property we can apply the continuity theorem~\ref{th.continuity}. We get for each $t \in [-t_0,t_0]$,
$$
\lim_{t \rightarrow +\infty} \rho(\leb,\t h \circ  \t \Phi_{s}^t)=  \rho(\leb,\t
h \circ \t \Phi^{t}).
$$
Moreover, from the compactly supported case we know that the morphism property holds for the compactly supported flow $(\Phi_{s}^t)_{t \in \RR}$:
$$
\rho(\leb,\t   h    \circ   \t \Phi_{s}^t)=\rho(\leb ,\t   h    )+   \rho(\leb,
\t \Phi_{s}^t).
$$
When $s$ tends towards $+\infty$ we get
$$
\rho(\leb,\t   h    \circ   \t \Phi^{t})=\rho(\leb ,\t   h    )+   \rho(\leb,
\t \Phi^{t})= \alpha +t . \rho(\leb,
\t \Phi^{1}).
$$
Thus in that subcase the morphism property also holds for the flow $(\Phi^{t})_{t \in \RR}$, and we conclude as in the compactly supported case.
\end{proof}

\paragraph{Remark}
\begin{enumerate}
\item
Note that the proof provides a periodic orbit whose rotation number is close to the rotation number $\alpha$ of the irrational pseudo-rotation $h$. However, in the case where the $(R_{\theta})$'s are the Euclidean rotations, we do not know if there must be a periodic orbit of rotation number $\alpha+\theta$ every time this number is rational. This is linked to the morphism property, see question~\ref{q.morphism}.
\item One sees on examples (that may even be conjugate to an irrational rotation)
that the two cases of the proof may actually occur.
\end{enumerate}

\appendix

\section{A hairy example}
\label{s.exemple-a-poils}
In this appendix we describe an example that shows why proposition~\ref{p.conjugue-a-translation}  is sharp.

\begin{figure}[htbp]
\label{f.hairy}
\begin{center}
\includegraphics[width=5cm]{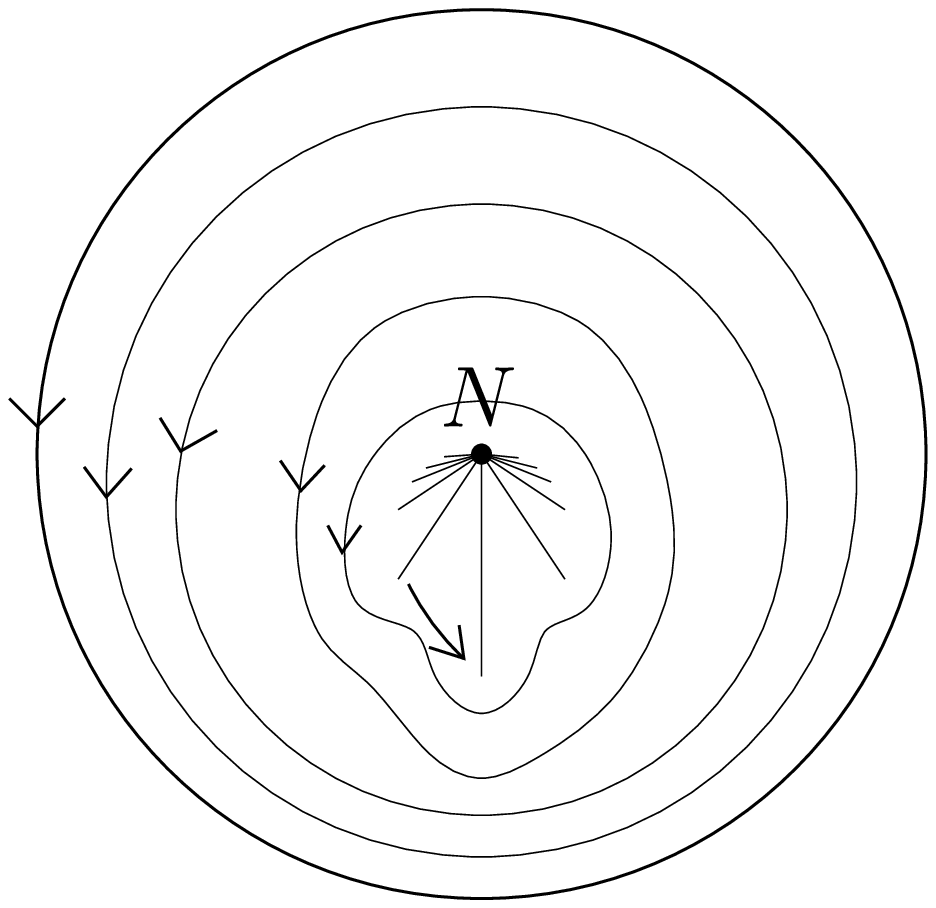}
\hspace{1cm}{\includegraphics[width=8cm]{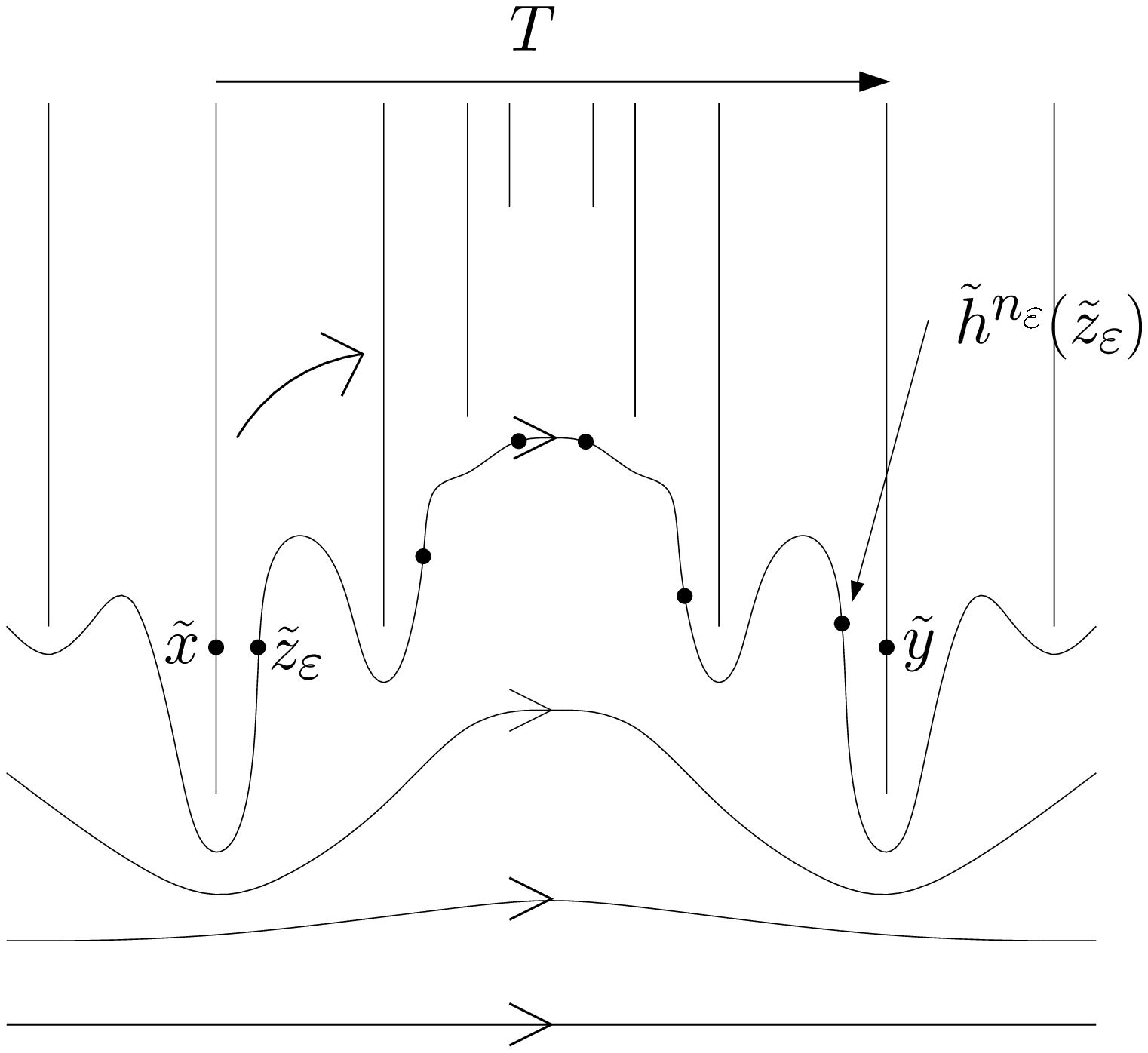}}
\caption{The hairy example}
\end{center}
\end{figure}

The left part of figure~\ref{f.hairy} shows the dynamics of a conservative homeomorphism $h$ of the disc.
The disc is foliated by circles, apart from the central set which is the union of a sequence $(I_{k})_{k \in \ZZ}$ of segments whose length tends to zero, having a single point $N$ in common. The homeomorphism fixes $N$ and sends $I_{k}$ to $I_{k+1}$.
Each circle is invariant, and the homeomorphism acts as a non-trivial rotation.
Note that the rotation number must tends to zero when the circles get closer to the hairy set.

Next we  see this disc as the upper half-sphere, and extend $h$ to the lower half-sphere by a rotation. Thus we get a conservative homeomorphism (again called $h$) of the sphere, with the two poles $N$ and $S$ as the only fixed points. We see $h$ as an element of $\homeo_\leb^+(\AA)$ with no fixed point. Then $h$ has a lift $\t h$ whose rotation set is equal to $]0,\alpha]$ for some $\alpha$. In particular, the points on the hairs are not recurrent, which explains why zero does not belong to the rotation set. Furthermore, one can prove that the homeomorphism $\t h$ is not conjugate to a translation. Indeed, consider a point  $\t x$  that projects on a hair, and $\t y = T(\t x)$ (see the right part of the figure). Then the couple $(\t x, \t y)$ is \emph{singular}: there exists points $\t z_{\varepsilon}$ arbitrarily near $\t x$, and arbitrarily large integers $n_{\varepsilon}$ such that $\t h^{n_{\varepsilon}}(\t z_{\varepsilon})$ is arbitrarily near $\t y$. This is a dynamical feature that distinguishes $\t h$ from a translation.

However there exists some essential topological lines $\gamma$  in $\AA$ that are disjoint from their image $h(\gamma)$ (take the projection of a vertical line in the right part of the figure). In general, when the rotation set is supposed to be included in the open interval $]0,1[$, we do not know if there always exist such a line. The line translation theorem~\ref{th.arc-translation} only works under the stronger assumption that the \emph{closure} of the rotation set is included in $]0,1[$.

As explained in paragraph~\ref{ss.existence-periodic}, we do not know either if a homeomorphism $h$ must have a fixed point when some lift $\t h $ has a rotation set equal to $[0,\alpha]$ for some positive $\alpha$.

\section{Conservative version of Schoenflies theorem}
\label{s.schoenflies}
\begin{theo}
\label{t.schoenflies}
Let $\alpha_{1}$, $\alpha_{2}$ be two simple arcs in the sphere $\SS^2$, and fix a homeomorphism $\phi$ from $\alpha_{1}$ onto $\alpha_{2}$. 
Then there exists an orientation-preserving homeomorphism $\Phi$ of $\SS^2$ which is an extension of $\phi$.
In addition to this, if the Lebesgue measure of $\alpha_{1}$, $\alpha_{2}$ is zero, then $\Phi$ can be chosen so that it preserves the Lebesgue measure.
\end{theo}
In this paper we use the theorem to straighten some essential topological lines $\Gamma$ of the annulus $\AA=\SS^1\times \RR$. However there is a slight difficulty coming from the fact that we have to deal with topological lines whose Lebesgue measure is not zero. We indicate here how to by-pass the problem.
The topological line $\Gamma$ comes with a homeomorphism $h$ of the annulus such that $\Gamma \cap h(\Gamma) = \emptyset$. We choose a neighbourhood $G$ of $\Gamma$ such that $G \cap h(G) = \emptyset$. Then one can find an essential  topological line $\Gamma'$, included in $G$, whose Lebesgue measure is zero
(for example, $\Gamma'$ can be piecewise affine). We replace $\Gamma$ by $\Gamma'$ before applying theorem~\ref{t.schoenflies}.

\begin{proof}[Idea of the proof]
We only indicate  how to get the second part from the first one.
Let $\Phi$ be an orientation-preserving homeomorphism that is an extension of $\phi$. 
Let $m$ denotes the image of the Lebesgue measure under $\Phi$. The measure $m$ has the following property: it is positive on each open set, it has no atom, it gives measure zero to the arc $\alpha_{2}$. According to a theorem of Oxtoby and Ulam (\cite{goff}), there exists an orientation-preserving  homeomorphism $\Psi$, that is the identity on $\alpha_{2}$, and that sends the measure $m$ on the Lebesgue measure
(actually, their theorem is stated on the square, but we can cut the sphere along the arc $\alpha_{2}$ to match this setting). The homeomorphism $\Psi \circ \Phi$ preserves the Lebesgue measure and is still an extension of $\phi$.
\end{proof}



\begin{thebibliography}{150}

\bibitem{anosov-katok}
Anosov, Dmitri and Katok, Anatole.
New examples in smooth ergodic theory. Ergodic diffeomorphisms.
\textit{Transactions of the Moscow Mathematical Society}
\textbf{23} (1970), 1--35.

\bibitem{beguin-leroux}
B\'eguin, Fran\c cois; Le~Roux, Fr\'ed\'eric.
Ensemble oscillant d'un  homéomorphisme de Brouwer, homéomorphismes de
Reeb. \textit{Bull. Soc. Math. France 131} \textbf{2} (2003), 149--210. 

\bibitem{beguin-crovisier-leroux-patou}
B\'eguin, Fran\c cois~; Crovisier, Sylvain~; Le~Roux, Fr\'ed\'eric and
Patou, Alice. 
Pseudo-rotations  of the  closed annulus:  variation on  a  theorem of
J. Kwapisz. Nonlinearity 17 (2004), no. 4, 1427--1453. 

\bibitem{birkhoff}
Birkhoff, George D.
\emph{Collected papers}, Vol II, Amer. Math. Soc., New York City, 1950. 

\bibitem{Cai}
Cairns, Stewart S.
An elementary proof of the Jordan-Schoenflies theorem.
\textit{Proc. Amer. Math. Soc.} \textbf{2} (1951), 860--867.


\bibitem{fathi-herman}
Fathi,  Albert  and  Herman,  Michael. Existence  de  difféomorphismes
minimaux.  Dynamical systems,  Vol.  I Warsaw.   \textit{Ast\'erisque}
\textbf{49} Soc. Math. France, Paris (1977), 37--59. 

\bibitem{fayad-katok}
Fayad, Bassam and Katok, Anatole. 
Constructions in elliptic dynamics.  \emph{Ergod. Th.  Dyn. Sys.} \textbf{24}  (2004),  no. 5, 1477--1520.

\bibitem{fayad-krikorian-vivier}
Fayad, Bassam~; Krikorian, Raphael and Vivier, Th\'er\`ese. 
In preparation. 

\bibitem{fayad-saprykina}
Fayad, Bassam and Saprykina, Maria.
Weak mixing diffeomorphisms on the disc and the annulus with arbitrary
Liouvillean rotation number on the boundary. Preprint. 


\bibitem{franks1}
Franks, John.
Generalizations of the Poincar\'e-Birkhoff theorem. 
\textit{Ann. of Math.} (2)  \textbf{128}  (1988),  no. 1, 139--151.

\bibitem{franks3}
Franks, John.
Rotation numbers for area preserving homeomorphisms of the open annulus.
\textit{Dynamical   systems  and   related  topics}   (Nagoya,  1990),
123--127, \textit{Adv. Ser. Dynam. Systems} 
\textbf{9}, World Sci. Publishing, River Edge, NJ (1991).

\bibitem{franks2}
Franks, John.
Area preserving homeomorphisms of open surfaces of genus zero.
\textit{New York J. Math.} \textbf{2} (1996), 1--19.

 \bibitem{goff} Goffman, Casper and Pedrick, George.
 A proof of the homeomorphism of Lebesgue-Stieltjes measure with Lebesgue measure.
 \textit{Proc. Amer. Math. Soc.} \textbf{52} (1975), 196--198. 

\bibitem{guillou1}
Guillou, Lucien.
Th\'eor\`eme de translation plane de Brouwer et g\'en\'eralisations du
th\'eor\`eme de Poincar\'e-Birkhoff. 
\textit{Topology} \textbf{33} (1994), no 2, 331--351.

\bibitem{guillou2}
Guillou, Lucien.
Unpublished text.


\bibitem{hall-turpin}
Hall, Glen and Turpin, Mark. 
Robustness of periodic point free maps of the annulus.
\textit{Topology Appl.} \textbf{69} (1996), no. 3, 211--215. 

\bibitem{handel}
Handel, Michael. 
A pathological area preserving $C\sp{\infty }$ diffeomorphism of the plane. 
\textit{Proc. Amer. Math. Soc.} \textbf{86} (1982), no. 1, 163--168.

\bibitem{herman1}
Herman, Michael.
Construction of some curious diffeomorphisms of the Riemann sphere.
\textit{J. London Math. Soc.} (2)  \textbf{34}  (1986),  no. 2, 375--384.

\bibitem{herman2}
Herman, Michael.
Some open problems in dynamical systems.
{\em  Proceedings of  the International  Congress  of Mathematicians},
Vol. II (Berlin, 1998). Doc. Math. 1998, Extra Vol. II, 797--808. 

\bibitem{kerekjarto}
Ker\'ekj\'art\'o, B\'ela. 
Sur     le    groupe     des    transformations     topologiques    du
plan.  \textit{Ann.  S.N.S.  Pisa},  t. \textbf{II},  Ser.  3  (1934),
393--400. 

\bibitem{kerekjarto2}
Ker\'ekj\'art\'o, B\'ela. 
\textit{Vorlesungen           \"uber           Topologie          (I),
Fl\"achentopologie}. Springer, Berlin, 1923. 


\bibitem{kwapisz}
Kwapisz, Jaroslaw.
A  priori degeneracy  of  one-dimensional rotation  sets for  periodic
point free  torus maps. \textit{Trans. Amer.  Math. Soc.} \textbf{354}
(2002), no 7, 2865--2895.  


\bibitem{lecalvez}
Le Calvez, Patrice.
Rotation numbers in the infinite annulus.
{\em Proc. Amer. Math. Soc.} \textbf{129} (2001), no. 11, 3221--3230.

\bibitem{lecalvez2}
Le Calvez, Patrice. 
Une version  feuillet\'ee équivariante  du théorème de  translation de
Brouwer. Preprint, 2004. 



\bibitem{leroux0}
Le Roux, Fr\'ed\'eric. 
\'Etude topologique de l'espace des hom\'eomorphismes de Brouwer. 
PhD thesis, Universit\'e Joseph Fourrier, Grenoble, 1997.

\bibitem{leroux}
Le Roux, Fr\'ed\'eric. 
Hom\'eomorphismes de surfaces - Th\'eor\`emes de la fleur de Leau-Fatou et
de la vari\'et\'e  stable. \textit{Ast\'erisque} \textbf{292} (2004). 

\bibitem{MisZie}
Misiurewicz, Michal and Ziemian, Krystyna.
Rotation sets for maps of tori.
\textit{J. London. Math Soc.} (2) \textbf{40} (1989), no.3, 490--506.

\bibitem{poincare}
Poincar\'e, Henri.
M\'emoire sur les courbes d\'efinies par une \'equation diff\'erentielle.
\textit{J. de Math Pures Appl.} S\'erie III \textbf{7} (1881), 375--422, \textbf{8} (1882), 251--296,
\textit{J. de Math Pures Appl.} S\'erie IV \textbf{1} (1885), 167--244, \textbf{2} (1886), 151--217.

 Also in \textit{{\OE}uvres de Henri Poincar\'e}, tome I, Gauthier Villars, Paris (1928), 3--44, 44--84,
90--158 and 167--222.


\bibitem{sauzet}
Sauzet, Alain.
\textit{Application des d\'ecompositions libres \`a l'\'etude des hom\'eomorphismes de surfaces.}
PhD thesis, Universit\'e Paris Nord, 2001.







\end{thebibliography}
\end{document}